\documentclass[11pt, reqno]{amsart}
\usepackage{graphicx,amsmath,amsfonts,amssymb,amsthm,color,framed,dsfont,hyperref,accents,tabularx,pgfplots}
\usepackage[numbers,sort]{natbib}
\pgfplotsset{compat=1.18}
\usepackage[T1]{fontenc} 
\usepackage[french,english]{babel}

\usepackage[margin=1in]{geometry}
\usetikzlibrary{external}
\usepackage{tikz}
\usetikzlibrary{arrows.meta, positioning}

\usepackage{physics}
\usepackage{bbm}
\usepackage[font={small,it}]{caption}

\newcommand{\R}{\mathbb{R}}

\newcommand{\Prob}{\mathbb{P}}

\newcommand{\E}{\mathbb{E}}

\newcommand{\Wass}{\mathbb{W}}
\newcommand{\Law}{\mathcal{L}}

\newcommand{\e}{\mathrm{e}}

\newcommand{\eps}{\varepsilon}
\newcommand{\ds}{\displaystyle}

\newcommand{\bnu}{\boldsymbol{\nu}}

\newcommand{\jl}{\eta}
\newcommand{\be}{\begin{equation}}
\newcommand{\ee}{\end{equation}}

\newtheorem{thm}{Theorem}[section]

\newtheorem{lemma}[thm]{Lemma}

\newtheorem{prop}[thm]{Proposition}
\newtheorem{assu}[thm]{Assumption}

\theoremstyle{remark}

\newtheorem{rem}[thm]{Note}
\AtBeginEnvironment{rem}{%
  \pushQED{\qed}%
}
\AtEndEnvironment{rem}{\popQED}
\numberwithin{equation}{section}

\title{Averaging for multi-scale McKean-Vlasov equation}

\author{Ashot Aleksian}
\author{Michela Ottobre} 
\address{Department of Mathematics, Heriot--Watt University, Edinburgh, EH14 4AS, United Kingdom} 

\email{a.aleksian@hw.ac.uk} \email{m.ottobre@hw.ac.uk}

\date{\today}

\begin{document}
\begin{abstract}
  We consider a system of Stochastic Differential Equations with two time scales, so called slow-fast system of SDEs.  As customary when studying slow-fast dynamics, we seek a coarse grained description of the system, obtained, in the limit of scale separation, by 'averaging out' the fast component of the dynamics (so-called averaging limit).  We study slow-fast systems where the fast variable is described by a process which is non-linear in the sense of McKean  i.e. by an SDE whose coefficients depend on the law of the process itself, while the slow variable is a linear SDE (linear in the sense of McKean). We show that, when the fast variable is a McKean-Vlasov process and, additionally, the coefficients of the SDE depend on the position of the slow process, the  limit of scale separation is not provided by the usual ansatz which one would expect by applying classical averaging paradigms. We identify the limit and provide complete proof. 

	\vspace{5pt}
    \noindent
    {\sc Keywords:}  Multiscale SDEs; McKean-Vlasov equations; averaging method. 
    
\vspace{5pt}
\noindent  
{\sc AMS Classification (MSC 2020):} 60H10; 60K35; 60J60; 37H05; 35J60. 
\end{abstract}

\maketitle

\section{Introduction}
This paper is concerned with the study of the averaging problem for two-scale, coupled SDEs, in the case in which the fast degrees of freedom are described by a McKean-Vlasov process.  More precisely, we consider the following system of SDEs, for the variables $X^\eps_t \in \R^d, $ $Y^\eps_t \in \R^m$: 
\begin{equation}
\label{eq:def:McKean-Vlasov_full}
\begin{cases}
    \dd{X_t^\varepsilon} = v(X_t^\varepsilon,Y_t^\varepsilon)\dd{t} + \sqrt{2 \kappa} \dd{B_t}, \\[6pt]
    \dd{Y_t^\varepsilon} = \dfrac{1}{\varepsilon} b(X_t^\varepsilon,Y_t^\varepsilon,\nu_t^\varepsilon)\dd{t}
    + \sqrt{\dfrac{2 \sigma}{\varepsilon}}\dd{W_t}, \\[8pt]
    \nu_t^\varepsilon = \mathcal L(Y_t^\varepsilon),
\end{cases}
\end{equation}
where $\nu_t^{\eps}$ is the law of $Y_t^{\eps}$,\footnote{Given a random variable $Z$, we denote by $\mathcal L(Z)$ the law of $Z$.}   $\kappa \geq 0$ is a non-negative constant,   
$(B_t)_{t\geq 0}$ and $(W_t)_{t\geq 0}$ are two independent Brownian motions in $\R^d$ and $\R^m$, respectively, and $0<\eps \ll 1$ is a parameter describing the scale separation between the slow variable $X_t^{\eps}$ and the fast variable $Y_t^{\eps}$. The initial conditions are taken to be deterministic: 
\[
    X_0^\eps = x_0 \in \R^d, \quad Y_0^\eps = y_0 \in \R^m \,, 
\]
and the drift functions are smooth functions with $v:\R^d\times \R^m \rightarrow \R^d$, $b:\R^d\times \R^m \times \mathcal P_2(\R^m) \rightarrow \R^m$, where $\mathcal P_2(\R^m)$ denotes the space of probability measures on $\R^m$ with finite second moment. \\
For the purposes of this paper, the important features  of  system \eqref{eq:def:McKean-Vlasov_full} are: the fast process is a McKean-Vlasov process (i.e. the evolution of $Y_t^{\eps}$ is non-linear in the sense that it depends on the law of the process $Y_t^{\eps}$ itself);\footnote{Throughout the paper, we use interchangeably the terms 'McKean-Vlasov process' and 'non-linear process in the sense of McKean'; and whenever, in short, we say non-linear process or non-linear SDE, we implicitly mean non-linear in the sense of McKean.} the drift coefficient $b$ depends on the position of the slow process;  and the slow process is described by an SDE (as opposed to an ODE), i.e. $\kappa>0$ is a strictly positive constant. To be more precise we treat both the case $\kappa =0$ and the case $\kappa>0$, see Section \ref{s:proof_main_result},  but the most relevant one to us is the latter.  As customary, we wish to identify the limit of the slow variable $X_t^{\eps}$ as $\eps$ goes to zero. This limit is routinely studied using the framework of stochastic averaging.  The purpose of this paper is to show that, in the setting of system \eqref{eq:def:McKean-Vlasov_full},  the limit which one would expect to obtain by using the usual averaging ansatz is not the correct one. In particular, if one was to use the classical averaging ansatz, then the limit process would be described by a linear SDE. We prove that, as long as $\kappa>0$, the limit dynamics is described by a McKean-Vlasov process. More precisely,  in Theorem~\ref{thm:main_thm} we  show that, as $\eps$ tends to zero, the slow process $(X^{\eps}_t)_{t\geq 0}$ converges to the process $(\bar X_t)_{t\geq 0}$, solution of the following McKean-Vlasov SDE:
\begin{equation}
    \label{eq:X_bar_rho_bar}
    \begin{cases}
        \dd{\bar{X}_t} = \bar{v}(\bar{X}_t, \bar{\rho}_t)\dd{t} + \sqrt{2\kappa} \dd{B_t}, \\
        \bar{\rho}_t = \Law(\bar{X}_t), \\
        \bar{X}_0 = x_0 \quad \text{a.s.},
    \end{cases}
    \end{equation}
    where the {\em averaged drift} $\bar{v}: \R^d \times \mathcal{P}(\R^d) \to \R^d$ is given by
    \begin{equation} \label{eqn: drift ansatz}
        \bar{v}(x, \rho) := \int_{\R^m} v(x, y) \Pi^{x, \rho}(\dd{y}) \,, 
    \end{equation}
and, for each $x\in\R^d$ and $\rho\in \mathcal P(\R^d)$,  $\Pi^{x, \rho}$ is is a probability measure on $\R^m$, found by solving equation  \eqref{eqn: Pi intro} below. We will discuss below and throughout the paper the precise meaning of the measure $\Pi^{x,\rho}$, and compare with the limit found when $\kappa=0$. For now we note that if $\kappa=0$ the dynamics \eqref{eq:X_bar_rho_bar} still makes sense (it is just the case that $\bar \rho_t$ will be a dirac); moreover,  we will show at the start of Section \ref{s:heuristics} that if $\kappa=0$ then \eqref{eq:X_bar_rho_bar}-\eqref{eqn: drift ansatz} reduces to the the averaged equation that one would expect (namely to \eqref{eqn:incorrect frozen process}-\eqref{eqn: limit process wrong} below). Hence  \eqref{eq:X_bar_rho_bar}-\eqref{eqn: drift ansatz} is truly only different from the expected averaging limit when $\kappa>0$. 

In order to better explain the result of this paper,  let us start by recalling the averaging paradigm, in the more standard context of linear SDEs. To this end, consider the system 
\begin{equation}
\label{eq:multiscale system intro}
\begin{cases}
    \dd{X_t^\varepsilon} = v(X_t^\varepsilon,Y_t^\varepsilon)\dd{t} + \sqrt{2 \kappa} \dd{B_t}, \\[6pt]
    \dd{Y_t^\varepsilon} = \dfrac{1}{\varepsilon} b(X_t^\varepsilon,Y_t^\varepsilon)\dd{t}
    + \sqrt{\dfrac{2}{\varepsilon} }\dd{W_t}, 
\end{cases}
\end{equation}
which (with some abuse of notation) can be viewed as a special case of \eqref{eq:def:McKean-Vlasov_full}, when the drift coefficient of the fast process does not depend on the law of the  process $Y_t^{\eps}$ itself (and we have fixed $\sigma=1$, which will be the case throughout the paper, as $\sigma$ does nos not play any role to our purposes). The averaging method can be intuitively understood by observing that,  as $\eps\rightarrow 0$, the process $Y^{\eps}$ will evolve much faster than the process $X^{\eps}$;\footnote{When this does not give rise to ambiguity, we use interchangeably the notations $(X_t)_{t\geq 0}$ and $X$ to denote a stochastic process.} in particular, the former will have reached equilibrium (assuming it has one) while the slow process has remained substantially unchanged. To formalise this idea one introduces an auxiliary process, $Y^{(\eps,x)}$, defined as the solution of the following SDE: 
\begin{equation}\label{eqn: auxil proc}
    \dd{Y_t^{(\eps,x)}} =\frac{1}{\eps} b(x, Y_t^{(\eps,x)}) \dd{t} + \sqrt{\frac{2}{\eps}} \dd{W_t} \, \,.
\end{equation}
Note that the above equation is obtained by "freezing" the value of the slow variable in the equation \eqref{eq:def:McKean-Vlasov_full}$_2$ for the fast process. 
Since $Y^{(\eps,x)}$ is decoupled from the evolution of the slow process, due to the time-scaling properties of Brownian motion, when looking at  process \eqref{eqn: auxil proc},  letting $t\rightarrow \infty$ for $\eps$ fixed is the same as letting $\eps\rightarrow 0$, for each $t$ fixed. Hence it is typical to fix the value of $\eps$ and consider 
the so-called {\em frozen process}, namely: 
\begin{equation}\label{eqn: frozen process intro}
    \dd Y_t^{(x)} = b(x, Y_t^{(x)}) \dd{t} + \sqrt{2} \dd{W_t} \, \,.
\end{equation}
Notice that the above is a family of processes parametrised by $x \in \R^d$,  and $x$ here is fixed deterministically.  

Assuming that, for each $x$ fixed, the process $Y^{(x)}$ admits a unique invariant measure $\pi^x$ -- which, we recall,  is the solution of the following partial differential equation: 
\be\label{eqn: pi^x classical case}
- \nabla_y \cdot \left( b(x,y) \pi^x\right) + \Delta_y \pi^x =0 \, ;
\ee  
and further more  assuming that $Y^{(x)}$ is ergodic with respect to such a measure, in the sense that
$$
\lim_{t\rightarrow \infty} \frac1t \int_0^t f(Y_s^{(x)}) \dd{s} = \int_{\R^m} f(y) \pi^x(\dd{y}), 
$$
for any function $f$ in a suitable class, one can prove that, as $\eps\rightarrow 0$, the process $X^{\eps}$ converges, at least on finite time horizons $[0,T]$, to the so-called {\em averaged dynamics} $\bar{\mathcal X} \subseteq \R^d$, solution of the SDE
\begin{equation}\label{eqn: standard limit intro}
   \dd{\bar{\mathcal X}_t} = \bar{v}(\bar{\mathcal X}_t) \dd{t} + \sqrt{2 \kappa} \dd{B_t} 
\end{equation}
where 
$$
\bar v(x) = \int_{\R^m} v(x,y) \pi^{x}(\dd{y})  \,. 
$$
Of course this result still holds when the system is fully coupled (at least assuming we are referring to weak convergence), i.e. when also the coefficients $\kappa$ and $\sigma$ depend on both the slow and the fast variables and in some cases can also be shown to hold on infinite time horizons \cite{crisan2026poisson}. 
Since the pioneering work of Khas'minskii  \cite{khas1968averaging},   the literature on averaging for stochastic dynamics has rapidly developed, \cite{FW, PS08, weinan2011principles, ethier2009markov},  starting from evolutions on compact spaces,  followed by the seminal papers \cite{pardoux2001poisson, pardoux2003poisson, pardoux2005poisson}, which paved the way for methods to prove the limit in non-compact settings. Here we mostly mention results on finite dimensional SDEs, which are  more closely relevant to this paper, but of course the literature on infinite dimensional multiscale systems is  fast developing, see e.g. \cite{cerrai2017averaging, hong2025averaging}. 

More closely related to what we study in this paper, recently various works have addressed both averaging and homogenization regimes for the case when either the slow process, the fast one, or both are McKean-Vlasov processes. In particular the works \cite{RSX21, WHY24, BS23, HLL26, xiang2026averaging} treat the case when the slow process is of McKean-Vlasov type (but not the fast one), with \cite{JZ25}  considering fractional Brownian motion as driving noise and \cite{HLL22} for infinite dimensional processes. The setting where the fast process is a McKean-Vlasov process has been studied \cite{SLLG25, XLLM21}.    One of the most general results we know of within this strand of literature can be found  in \cite{LietalPoisson}, where the drift of the slow process depends on $X_t^{\eps}, Y_t^{\eps}, \mathcal L (X_t^{\eps}), \mathcal L (Y_t^{\eps})$, while the drift of the fast process depends on $Y_t^{\eps}, \mathcal L (X_t^{\eps})$ and $ \mathcal L (Y_t^{\eps})$, but it does not depend on $X_t^{\eps}$ itself. Recently, the work \cite{XLLM21} has tried to find the averaging limit when the drift of the fast process depends on $X_t^{\eps}$. We believe their result to be incorrect (at least in the case $\kappa>0$) and correct it here.  More comments on this in Section \ref{s:heuristics} and in Note \ref{note:explnation error} we explain in detail why we think the ansatz for the limit in \cite{XLLM21} is incorrect. 

Coming back to our system of interest, system \eqref{eq:def:McKean-Vlasov_full}, if we were to follow brute force the `recipe' that we have highlighted above, then we would be tempted to consider the frozen process
\begin{equation}\label{eqn:incorrect frozen process} 
    \dd{Y_t^{(x)}} = b(x, Y_t^{(x)}, \nu_t^{(x)}) \dd{t} + \sqrt{2}  \dd{W_t} \, , \quad \nu_t^{(x)} = \mathcal L\left(Y_t^{(x)}\right)
\end{equation}
i.e. the process obtained by  freezing the value of the slow process in the equation for the fast process. 
The (density of the) stationary measure associated with the above process (if it exists) is the solution of the PDE
\be\label{eqn: limit measure wrong}
- \nabla_y \cdot \left[ b(x,y,\nu^{(x)}) \nu^{(x)}\right] + \Delta_y \nu^{(x)} =0 \,.
\ee
Based on this reasoning one would conclude that the averaged dynamics is given by 
\be\label{eqn: limit process wrong}
    \dd{\bar{\mathcal X}_t} 
    = \bar v(\bar{\mathcal X}_t) \dd{t} 
    + \sqrt{2 \kappa }\,  \dd{B_t} \, , \quad 
    {\mbox{ with }} 
    \bar{v}(x) = \int_{\R^m} v(x,y) \nu^{(x)}(\dd{y}) \,.
\ee
In particular the averaged dynamics would be linear in the sense of McKean. 
We claim that the measure $\nu^{(x)}$ which we have defined in the above is not the measure that one should average against in order to obtain the correct coarse grained description. In particular, because of the nonlinear nature of the fast process,  to obtain a process which plays a role similar to the one played by the frozen process, we cannot just freeze the value of $X_t^{\eps}$ in the equation for the fast variable. In the next section we will show some heuristics (and then, in subsequent sections, complete  proofs) to provide intuition about the fact that, since the fast process is a McKean-Vlasov process (and the slow variable is subject to its own source of noise), when building the analogue of the frozen process, we need to consider a family of processes which are not just parametrised by the possible outcomes $x$ of the slow process; indeed, we will show that both the outcome and  the law of the random variable we `freeze' will have an effect on the fast dynamics (see next section). It is intuitively for this reason that, when building the analougue of the frozen process,  we need to consider a family of processes parametrised by an entire random variable. More explicitly, we will consider  processes of the form 
\begin{equation} \label{eqn: our frozin proces intro}
    \dd{\mathcal Y_t^{\xi}} = b(\xi, \mathcal Y_t^{\xi}, \boldsymbol{\nu}_t^{\xi}) \dd{t} 
    + \sqrt{2} \dd{\widetilde{W}_t}\, , \quad
    \bnu_t^{\xi} = \mathcal L(\mathcal Y_t^{\xi}) \,, 
\end{equation}
where in the above $\xi$ is a random variable in $L^2(\Omega; \R^d)$ and $\widetilde{W}$ is a standard Brownian motion.\footnote{The precise definition of this process is given in   \eqref{eq:frozen_process_1_gen}, where it is denoted by $\mathcal Y^{\xi, \zeta}$. Here we omit the index $\zeta$, as for the purposes of this introduction we focus on the dependence on $\xi$ only.}  Hence the above is a family of processes parametrised by the random variable $\xi$. Note that while $\boldsymbol{\nu}^{\xi}$ defined in the above depends on the law of $\xi$ (in the context of the overall proof, where the role of $\xi$ is played by the variable $\bar X_t$ at  certain fixed times, $\boldsymbol{\nu}^{\xi}$ will depend on the law of $\bar X_t$), when we consider $\nu^{(x)}$ in \eqref{eqn:incorrect frozen process} $x$ is fixed, so there is no such a dependence. With this in mind, let $\rho$ denote the law of $\xi$. Then the joint law    $\boldsymbol{\mu}_t^{\xi} = \mathcal L (\xi, \mathcal Y_t^{\xi})$  
 can be written as 
\be\label{eqn: mu intro}
    \boldsymbol{\mu}_t^\xi (\dd{x}, \dd{y}) = \rho(\dd{x}) \varsigma_t^{x, \rho}(\dd{y}), 
    \qquad 
    \varsigma_t^{x, \rho}:= \mathcal L \left(\mathcal Y_t^{\xi} \vert \xi=x \right) \,,  \footnote{Most of the measures used in this paper will have densities with respect to the Lebesgues measure. When this is the case, with abuse of notation, we denote by the same letter both the measure and the density, i.e. we will often write $\nu(\dd y)=\nu(y) \dd y \,.$  }
\ee
where $\varsigma_t^{x, \rho}$ is defined for $\rho$-almost every $x$. For $\rho-$ almost every $x$, the conditional law $\varsigma_t^{x, \rho}$  satisfies the following  PDE (see Section \ref{subsec: proof of prop 3.5} and Appendix \ref{appendix: conditional PDE}), 
\be\label{eqn: evol conditional law intro}
\partial_t \varsigma_t^{x, \rho}(y) = - \nabla_y \cdot \left[ b(x,y,\bnu_t^{\xi}) \varsigma_t^{x, \rho}(y)
\right] + \Delta_y \, \varsigma_t^{x, \rho}(y) \,,
\ee 
where again $\bnu^{\xi}$ is the marginal $\ds \bnu_t^{\xi}(y) = \int_{\R^d} \boldsymbol{\mu}_t^\xi({x}, y)\dd x = \int_{\R^d} \varsigma_t^{z, \rho}(y)  \rho(z) \dd{z}$. The stationary solution $\Pi^{x, \rho}$ of the above PDE (which exists and is unique by the assumptions we will make, see Section \ref{s:Assumptions_main_result}) is then the solution to
\be\label{eqn: Pi intro}
- \nabla_y \cdot \left[ b\left( x,y, \int_{\R^d} \Pi^{z, \rho}(y)  \rho(z) \dd{z} \right) \Pi^{x, \rho}(y)
\right] + \Delta_y \, \Pi^{x, \rho}(y) = 0 \,.
\ee
This motivates using $\Pi^{x, \rho}$ in the averaged equation \eqref{eq:X_bar_rho_bar}. 
Let us emphasize that solving  \eqref{eqn: Pi intro}  means finding, for each fixed measure $\rho \in \mathcal P_2(\R^d)$, a family $(\Pi^{x, \rho})_{x \in \R^d} \subset \mathcal P_2(\R^m)$ which solves the above. The same can be said about \eqref{eqn: evol conditional law intro}, the solution of which is a family indexed by $x \in \R^d$.  To prove well-posedness of \eqref{eqn: Pi intro} we will recast it as a system, and solve it via a fixed point argument, see Proposition \ref{prop:Pi_x_rho_existence}.\\

The methods we employ to prove the convergence of $X^{\eps}$ to $\bar X$ are based on a classic probabilistic approach, which will require proving an ergodic result for processes of the form \eqref{eqn: our frozin proces intro}, see Proposition \ref{prop:ergodicity_of_Y}.  
We note in passing that as a consequence of our results,  for $\rho$-a.e. $x$,  $\varsigma_t^{x,\rho}$ converges to $\Pi^{x,\rho}$ as $t \rightarrow \infty$, see Note \ref{note: convergence}.   

Finally, in \eqref{eq:def:McKean-Vlasov_full} we chose to restrict to constant diffusion coefficients, but generalization are possible (with the most straightforward being the case of $\kappa$ and $\sigma$ dependent of the slow variable $X_t^{\eps}$ and elliptic). Similarly, very little would change in the proof if both the drifts $b$ and $v$ depended on the law of the slow variable. We don't consider these cases here for brevity and because we wish to convey the bulk of the novelty of the ansatz \eqref{eq:X_bar_rho_bar}-\eqref{eqn: Pi intro} in a clear manner, without unnecessary technicalities obscuring the main point.

Throughout this paper, we work on a probability space 
\( (\Omega,\mathcal F,\mathbb P) \) 
that is rich enough to support all random variables and stochastic processes appearing in the sequel.

The paper is organised as follows. In the next Section \ref{s:heuristics} we give some heuristics to explain the main result of this paper and provide more comparison to literature. Section \ref{s:Assumptions_main_result} contains the statement of assumptions and of main results. Section \ref{s: proofs auxiliary results} and Section \ref{s:proof_main_result} are devoted to proofs. Appendix \ref{appendix A} contains miscellaneous proofs and a further heuristic argument, based on looking at the behaviour of the joint law $\eta^{\eps}$ of the process \eqref{eq:def:McKean-Vlasov_full}.

\section{Heuristics}
\label{s:heuristics}
In the introduction we have highlighted that the reason why our ansatz \eqref{eq:X_bar_rho_bar} for the coarse grained dynamics is not the expected one depends on a combination of three facts:   the fast process $Y_t^{\eps}$ is of McKean-Vlasov type; it depends on the position of the slow process $X_t^{\eps}$  (i.e.  the drift coefficient $b$ in \eqref{eq:def:McKean-Vlasov_full}  depends on $X_t^{\eps}$); the slow process is subject to its own source of noise, i.e. $\kappa>0$. 

So, let us start by formally showing how our result reduces to known cases when either one of these three assumptions does not hold.  
 Suppose the fast process $Y^{\eps}$ is not of McKean-Vlasov type, i.e. the drift coefficient $b$ in \eqref{eq:def:McKean-Vlasov_full} does not depend on $\nu_t^{\eps}$;  then it is easy to see that equation \eqref{eqn: Pi intro} reduces precisely to equation \eqref{eqn: pi^x classical case} and so, as expected, the coarse grained dynamics $\bar X$ would coincide with the dynamics $\bar{\mathcal X}$ in \eqref{eqn: standard limit intro}.

Now consider the case when  the fast process $Y^{\eps}$ is of McKean-Vlasov type but does not depend on $X^{\eps}$, i.e. when the drift coefficient $b$ in \eqref{eq:def:McKean-Vlasov_full} does not depend on $X_t^{\eps}$. Then the fast process becomes decoupled of the slow process and it is well known (for example as a corollary of the results of \cite{LietalPoisson}) that the limit coarse grained dynamics is obtained by averaging the drift against the invariant measure of the fast process (which now is not parametrised in any way by the slow process). To see that our result is consistent with this, note first that our frozen process in  \eqref{eqn: our frozin proces intro} would not depend on $\xi$ and would reduce to the process 
\begin{equation} \label{eqn: simpler frozen}
    \dd{\mathcal Y_t} = b( \mathcal Y_t, \nu_t) \dd{t} 
    + \sqrt{2}  \dd{W_t}\, , \quad
    \nu_t = \mathcal L(\mathcal Y_t) \,. 
\end{equation}
Since the drift in the above is not random, i.e. there is no $\xi$, the measure $\boldsymbol{\mu}_t$ in \eqref{eqn: mu intro} would reduce to a measure $\mu_t$ on $\R^m$ only and in particular it reduces to simply $\mu_t= \mathcal L(\mathcal Y_t)= \nu_t$ and there is no conditional law to be considered. As $t\rightarrow \infty$, such a measure would precisely converge to the invariant measure of the fast process. Differently put, looking at equation \eqref{eqn: Pi intro}, if there is no $\rho$ and no $x$, such an equation coincides precisely with the equation for the stationary measure of \eqref{eqn: simpler frozen}. 

Finally, suppose the fast process is as in \eqref{eq:def:McKean-Vlasov_full} but the diffusion coefficient of the slow process vanishes, i.e. $\kappa=0$. In this case  the averaged dynamics boils down to the one given by \eqref{eqn: limit process wrong}, with the averaged drift resulting from averaging with respect to the measure $\nu^{(x)}$ in \eqref{eqn: limit measure wrong}. To see  this, notice that when $\kappa=0$ in \eqref{eq:X_bar_rho_bar}  we  have  $\bar \rho_t = \delta_{\bar X_t}$, as $\bar X$ would be deterministic. Hence the averaged coefficient $\bar v$ in \eqref{eqn: drift ansatz} is known once the family of measures $\Pi^{x, \delta_x}$ is known. It is easy to see that when, in \eqref{eqn: Pi intro}, we take $\rho$ to be a dirac delta centered in $x$, the resulting measure satisfies the same equation as $\nu^{(x)}$, equation \eqref{eqn: limit measure wrong},  and hence coincides with $\nu^{(x)}$. 

\smallskip
We now move on and present some heuristic arguments to give intuition about the main result. We provide two heuristics: the first, in Subsection \ref{ss:time-discretisation_particle_system}, is mostly aimed at giving intuition about why, in our setting,  one would expect the coarse grained dynamics to be a non-linear process despite the fact that the slow variable in \eqref{eq:def:McKean-Vlasov_full} is (seemingly) a linear process in the sense of McKean; the second, in Subsection \ref{ss: second heuristics}, provides an argument which is closer to the one followed in the proof.

\subsection{A first heuristic argument: particle approximation of \eqref{eq:def:McKean-Vlasov_full}}
\label{ss:time-discretisation_particle_system}

In order to gain some intuition on why the coarse grained dynamics \eqref{eq:X_bar_rho_bar} is a non-linear SDE (despite the fact that the slow process $X^{\varepsilon}$ in \eqref{eq:def:McKean-Vlasov_full}$_1$ is, seemingly, a linear process) we  consider the  particle approximation of system \eqref{eq:def:McKean-Vlasov_full}. To write such a particle approximation,  it is important to observe that  the pair $(X^\eps, Y^\eps)$  can be regarded as a standard McKean-Vlasov process in $\R^{d \times m}$; that is, system~\eqref{eq:def:McKean-Vlasov_full} can be written as:
\begin{equation*}
    \begin{pmatrix}
        \dd{X_t^\eps} \\
        \dd{Y_t^\eps}
    \end{pmatrix}
    = 
    \begin{pmatrix}
        1 & 0 \\
        0 & 1/\eps
    \end{pmatrix} 
    V \big(X_t^\eps, Y_t^\eps, \jl^\eps_t \big) \dd{t} + 
    \begin{pmatrix}
        \sqrt{2} & 0 \\
        0 & \sqrt{2/\eps}
    \end{pmatrix}    
    \begin{pmatrix}
        \dd{B_t} \\
        \dd{W_t} 
    \end{pmatrix}  
\end{equation*}
where $\jl_t^\eps$ denotes the joint law $\Law(X_t^\eps, Y_t^\eps)$ and  $V$ is a function of $(x, y, \jl) \in \R^d \times \R^m \times \mathcal{P}_2(\R^d \times \R^m)$, namely
\begin{equation*}
    V(x, y, \jl) =     
    \begin{pmatrix}
        v(x, y) \\
        b\left(x, y, \int_{\R^d}\jl(x,y) dx\right) 
    \end{pmatrix}.
\end{equation*}
With this in mind, the particle system approximation for the pair $(X^\eps, Y^\eps)$ is given by
\begin{equation}
\label{eq:sys_of_particles_approx}
\begin{cases}
    \dd{X_t^i} = v\big(X_t^i,Y_t^{i}) \dd{t} + \sqrt{2}\dd{B_t^i}, \quad i \in \{1, \dots, N\},\\[0.3em]
    \dd{Y_t}^{i} = \dfrac{1}{\varepsilon}
    b\big(X_t^i, Y_t^{i},\bar{\nu}_t^N \big) \dd{t}
    + \sqrt{\dfrac{2}{\varepsilon}}\,\dd{W_t^i}, \quad i \in \{1, \dots, N\},\\[0.6em]
    \bar{\nu}_t^N = \frac{1}{N} \sum_{i = 1}^N \delta_{Y^{i}_t}.
\end{cases}
\end{equation}
As expected from propagation of chaos, each pair \((X^i,Y^i)\) (here  we view each pair \((X^i,Y^i)\) as being one particle) becomes asymptotically independent of the remaining pairs as \(N \to \infty\) and converges to the solution of \eqref{eq:def:McKean-Vlasov_full}. We don't prove this fact, which follows by standard theory, see e.g. \cite{chaintron2022propagation}. We wish to use this particle system approximation of \eqref{eq:def:McKean-Vlasov_full} to explain why the multiscale limit of \eqref{eq:def:McKean-Vlasov_full} is a McKean-Vlasov process rather than a linear  SDE.

System \eqref{eq:sys_of_particles_approx} consists of $N$ particles that evolve in $\R^{d \times m}$ and the  second component of each particle moves very fast. It is a linear system (in the sense of McKean)  and can therefore in principle be approximated using a standard Euler scheme. However if one was to use a Euler scheme then, as is well known, \cite[Chapter 9]{weinan2011principles}, in order to capture the behaviour of the fast dynamics one would need to choose a time-step which is too small (order $\eps$ or smaller), with the result that the algorithm typically gets stuck.

\begin{figure}[b]
\centering

\includegraphics{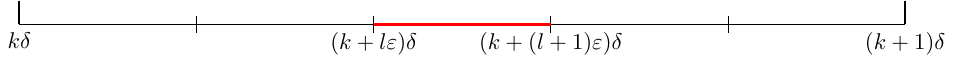}

\caption{A single step of size $\delta$ on the slower time scale corresponds to approximately $1/\varepsilon$ steps of size $\varepsilon \delta$ on the faster time scale.}
\label{fig:time_scales}
\end{figure}

To circumvent this problem one often takes the following approach \cite{weinan2011principles}.  Divide the time-interval $[0,T]$ in time steps of size $\delta$, which will be a small time step on the time-scale of the slow variable. Then the position of the $i$-th  slow particle is approximated by $X^i_{k \delta}$ on the entire interval $[k \delta , (k + 1) \delta )$.   Assuming  that $1/\eps$ is an integer, for each $i \in \{1, \dots, N\}$ the sequence $\{X^i_{k\delta}\}_{k\geq 0}$ is produced as follows: given $X^i_{k\delta}$ we evolve the fast process, i.e. we calculate 
 $$
Y^{i}_{(k + (l + 1) \eps) \delta } = Y^{i}_{(k + l \eps )\delta } +
    \frac{1}{\eps}b\big(X^i_{k \delta }, Y^{i}_{(k + l \eps) \delta }, \bar{\nu}^N_{(k + l \eps) \delta } \big) \delta
    + \sqrt{\frac{2 \delta}{\eps} } \alpha^i \,, 
 $$
 where $\ell = 0,1, \dots, \eps^{-1}$ and the $\alpha^i$'s are i.i.d. standard Gaussian random variables.    
Then the position of the slow process is updated: 
\be\label{eqn:slow update}
    \displaystyle X^i_{(k + 1) \delta} = X^i_{k \delta } + \left[ \eps \sum_{l = 1}^{1/\eps} v\big(X^i_{k \delta }, Y^{i}_{(k + l \eps) \delta }) \right]\delta  + \sqrt{2 \delta } \widetilde{\alpha}^i \,,  
\ee  
again with the $\widetilde{\alpha}^i$'s being i.i.d. standard Gaussian random variables, independent of the $\alpha^i$'s. 
This procedure provides a way to separate the two time scales: first we update the positions of $Y^i$'s $1/\eps$ times over a time interval of size $\delta$ and then we update the $X^i$'s using the cumulative effect of all these small changes. 

Now we can intuitively see how, in general, the dynamics of \(X^i\) is influenced by the other \(X^j\)'s -- which, in the continuum limit, would amount to a dependence on the law of the slow variable.  Indeed, looking at \eqref{eqn:slow update}, the evolution of $X^i$ depends on the evolution of $Y^i$; the latter depends on the empirical measure \(\bar{\nu}_t^N\), which contains the positions of all other \(Y^j\)'s. Since these variables are themselves driven by the corresponding \(X^j\)'s, the evolution of \(Y^i\) indirectly incorporates information about the entire particle system. This mechanism is precisely what gives rise to the law-dependence in the effective equation obtained in the averaging limit. Of course this point of view only provides intuition (and neglects the issue of the $N\rightarrow \infty$ limit and the $\eps\rightarrow 0$ limit potentially not commuting). 

Observe that this phenomenon disappears in some important special cases. First, if the dynamics of \(Y\) does not depend on its own law then the corresponding approximating particle \(Y^i\) evolves independently of the other particles. Second, if the drift $b$ of the fast process  does not depend on the slow variable \(X\), then the dynamics of the \(Y^j\)'s carries no information about the positions of the \(X^j\)'s. Although the \(Y^j\)'s may still interact through the empirical measure \(\bar{\nu}^N\), this interaction is independent of the slow variables. In either case,  the effective equation is no longer expected to exhibit the dependence on its own law, i.e. it is no longer expected to be of McKean-Vlasov type.

We remark in passing that  one may think that the particle approximation to \eqref{eq:def:McKean-Vlasov_full} should look like:
\begin{equation}
\label{veq:wrong_sys_of_particles}
\begin{cases}
    \dd{X_t} = v\big(X_t,Y_t^{1}) \dd{t} + \sqrt{2}\dd{B_t},\\[0.3em]
    \dd{Y_t}^{i} = \dfrac{1}{\varepsilon}
    b\big(X_t, Y_t^{i},\bar{\nu}_t^N \big) \dd{t}
    + \sqrt{\dfrac{2}{\varepsilon}}\,\dd{W_t^i}, \text{ for any } i \in \{1, \dots, N\}, \\[0.6em]
    \bar{\nu}_t^N = \frac{1}{N} \sum_{i = 1}^N \delta_{Y^{i}_t}.
\end{cases}
\end{equation}
Yet, this is not correct. Indeed, we know that such a system of particles converges to a conditional McKean-Vlasov process of the form 
\begin{equation*}
\begin{cases}
    \dd{X_t} = v(X_t,Y_t)\dd{t} + \sqrt{2} \dd{B_t}, \\[6pt]
    \dd{Y_t} = \dfrac{1}{\varepsilon} b(X_t,Y_t, \widetilde{\nu}_t)\dd{t}
    + \sqrt{\dfrac{2}{\varepsilon}}\dd{W_t}, \\[8pt]
    \widetilde{\nu}_t = \mathcal L(Y_t| \mathcal{F}^B_t),
\end{cases}
\end{equation*}
and not to \eqref{eq:def:McKean-Vlasov_full} (see e.g. \cite{CZ16}).  Here, $\mathcal{F}^B$ denotes the filtration with respect to the Brownian motion $B$ that drives the slow process $X$. Intuitively, this is due to the fact that the process $X$ in \eqref{veq:wrong_sys_of_particles} acts on the particles $Y^i$ as a common noise so in the limit we obtain conditional propagation of chaos.

\subsection{Second heuristic argument.}\label{ss: second heuristics}

In this section, we heuristically derive the form of the effective dynamics.   
In other words, assuming the slow process $X^{\eps}$ in \eqref{eq:def:McKean-Vlasov_full} converges, as $\eps \rightarrow 0$, to some process $\bar X$, we want to give a formal argument to show that $\bar X$ should coincide with the solution of \eqref{eq:X_bar_rho_bar}. While informal, the argument we illustrate here is substantially the backbone of the proof of our main result Theorem \ref{thm:main_thm}.  

As before, let us call $\delta$ a small time step on the slow time scale. 
From \eqref{eq:def:McKean-Vlasov_full}, for any  time $t$ in an interval  $[k \delta , (k + 1) \delta )$, we can write 
\begin{align}
    X_t^{\eps} & = X_{k\delta}^{\eps}+ \int_{k\delta}^t v(X_s^{\eps}, Y_s^{\eps}) ds + \sqrt{2 \kappa} \, dB_t \nonumber\\
    & \simeq X_{k\delta}^{\eps}+ \int_{k\delta}^t v(X_{k\delta}^{\eps}, Y_s^{\eps}) ds + \sqrt{2 \kappa} \, dB_t  \label{use1}
\end{align}
As for the fast process $Y_s^{\eps}$, we will show in Section \ref{s:proof_main_result} that it is well approximated by a process $\widetilde Y_t^{\delta, \eps}$ defined, for $t \in [k\delta, (k+1)\delta )$, as 
$$
\widetilde Y_t^{\delta, \eps} = \widetilde Y_{k\delta}^{\delta, \eps} + \frac{1}{\eps} \int_{k\delta}^t b(\bar X_{k\delta}, \widetilde Y_s^{\delta, \eps}, \widetilde{\nu}_s^{\delta ,\eps})  \, ds + \sqrt{\frac{2}{\eps}}  dW_t , \qquad  \widetilde{\nu}_t^{\delta ,\eps} = \mathcal L(\widetilde{Y}_t^{\delta ,\eps}) \,.
$$
Note that the above process is nonlinear in the sense of McKean, and its drift is random, because of the presence of the random variable $\bar X_{k\delta}$. The law of $\widetilde Y_t^{\delta, \eps}$ will then depend on the law (as well as on the simple outcomes) of the random variable $\bar X_{k\delta}$.  In other words, the process $\widetilde{Y}^{\delta, \eps}$ has the same structure of the process \eqref{eqn: our frozin proces intro}. Moreover,  with this definition, $\widetilde Y^{\delta, \eps}_t = \mathcal{Y}^{\xi}_{(t - k\delta )/\eps}$ almost surely for $t \in [k \delta , (k + 1)\delta )$, assuming $\mathcal Y^{\xi}_0=  Y^{\eps}_{k\delta}$ and taking $\xi$ to be $\bar X_{k\delta}$  (and modulo other caveats and subtleties which we will clarify in Section \ref{s:proof_main_result}); that is $\widetilde{Y}^{\delta, \eps}$ is substantially only a time-change with respect to $\mathcal Y^{\xi}$.    {\it Studying the behavior of $\widetilde{Y}^{\delta, \eps}$ in $[k \delta, (k + 1) \delta)$ in the regime $\eps \ll 1$ and for $t$ fixed, amounts to studying the long-term behavior of processes $\mathcal{Y}^{\xi}$ of the form \eqref{eqn: our frozin proces intro} }, which motivates considering  processes of the type \eqref{eqn: our frozin proces intro} as a frozen process. 

We emphasize again that the the  big difference between the class of frozen process we consider and say \eqref{eqn: frozen process intro} or \eqref{eqn:incorrect frozen process} is that, in our case, the frozen process is a McKean-Vlasov process with random drift term and therefore it depends on the law of $\xi$ in a non-trivial way. 
This is crucial for the evolution of the law of $\mathcal{Y}^{\xi}$: it has to take into account the distribution of $\xi$. 

Going back to \eqref{use1}, with the above in mind and recalling that we assume $X_t^{\eps} \simeq \bar X_t$,  we can then write 
\begin{align}
    X_t^{\eps} & \simeq  X_{k\delta}^{\eps}+ \int_{k\delta}^t v(X_{k\delta}^{\eps}, \widetilde{Y}_s^{\delta, \eps}) ds + \sqrt{2 \kappa} \, dB_t  \nonumber \\
    & = X_{k\delta}^{\eps}+ \int_{k\delta}^t v(X_{k\delta}^{\eps}, \mathcal{Y}_{(s-k\delta)/\eps}^{\bar X_{k\delta}}) ds + \sqrt{2 \kappa} \, dB_t \nonumber\\
    & = X_{k\delta}^{\eps}+ \eps \int_{0}^{(t-k\delta)/\eps} v(X_{k\delta}^{\eps}, \mathcal{Y}_{s}^{\bar X_{k\delta}}) ds + \sqrt{2 \kappa} \, dB_t \nonumber\\
    & \simeq \bar X_{k\delta}+ \eps \int_{0}^{(t-k\delta)/\eps} v(\bar X_{k\delta}, \mathcal{Y}_{s}^{\bar X_{k\delta}}) ds + \sqrt{2 \kappa} \, dB_t
    \label{use2}
\end{align}

From the last equation, since $\delta$ is fixed and we want to let $\eps$ to zero, it is clear that we need to study  limits of the form
$$
\lim_{t\rightarrow \infty}\frac{1}{t} \int_0^t f(\xi, \mathcal Y^{\xi}_s)  ds \,.
$$
It should be noted that the limit of the above time averages is a random quantity. We will show in Section \ref{s: proofs auxiliary results} that such a limit can be studied by considering the conditional law  $\varsigma^{x, \rho}_t = \Law(\mathcal{Y}^{\xi, \zeta}_t | \xi = x)$; an important step is then to show that $\varsigma^{x, \rho}_t$   converges to $\Pi^{x, \rho}$, defined in \eqref{eqn: Pi intro}, where we recall $\rho = \Law(\xi)$.  
 We will prove  (see Section~\ref{s:Assumptions_main_result} Proposition~\ref{prop:ergodicity_of_Y}), that for any regular enough function $f$ the following holds 
\begin{equation}
\label{eq:aux:ergodic_result}
    \frac{1}{t} \int_0^t f(\xi, \mathcal{Y}^{\xi}_s) \dd{s} \xrightarrow[t \to \infty]{} \int_{\R^m} f(\xi, y) \Pi^{\xi, \rho}(\dd{y}),    
\end{equation}
almost surely and in $L^2 (\Omega; \R)$.  It is now clear that, taking the limit $\eps \to 0$, and using the ergodic result \eqref{eq:aux:ergodic_result}, we obtain

\[
    \bar X_{t} \simeq X^{\eps}_t  \simeq \bar X_{k\delta}+ \delta \int_{\R^m} v(\bar X_{k \delta}, y) \Pi^{\bar X_{k \delta }, \Law(\bar X_{k \delta })} (\dd{y})  ,
\]
which is precisely  \eqref{eq:X_bar_rho_bar} (for $t \in [k\delta, (k+1)\delta)$).

\begin{rem}
\label{rem:Y_xi_zeta_parameters}
We emphasize again  that $\xi \in L^2(\Omega; \R^d)$  should be viewed as a \emph{parameter} of the process $\mathcal{Y}^{\xi}$. In contrast to the case of a classical linear SDEs, one cannot construct a family of processes
\(
z\mapsto \mathcal{Y}^{z}
\)
and then recover $\mathcal{Y}^{\xi}$ by substituting outcomes of the random variables $\xi$, that is, by considering the process
\(
\omega \mapsto \mathcal{Y}^{\xi(\omega)}(\omega).
\)
Such a representation is incorrect, since the dynamics are non-linear: the process depends on its own law. Consequently, $\xi$  cannot be interpreted as random parameter in the usual pathwise (outcome-wise) sense. 
Rather, it serves only as a parameter to define the process \eqref{eqn: our frozin proces intro}, and the dependence on the elementary outcome $\omega$ arises solely through
\(
\omega \mapsto \mathcal{Y}^{\xi}(\omega).
\) A different way of expressing this fact is by noting that the process $\mathcal Y^{\xi}_t \vert \xi =x$ is not the the same as the process $\mathcal Y^{(x)}$  in \eqref{eqn:incorrect frozen process} but rather the process $\mathcal{Y}^{\xi\vert x}$ which solves the following SDE
\begin{equation}\label{eqn:conditional process}
   \dd \mathcal{Y}^{\xi\vert x}_t = b(x, \mathcal{Y}^{\xi\vert x}_t, \bnu_t^{\xi}) dt + \sqrt 2 \,  d\widetilde{W}_t \,. 
\end{equation}
where we emphasize that in the above the third argument of $b$ is the second marginal $\bnu^{\xi}_t$ of the joint law $\boldsymbol{\mu}_t^{\xi}$ of $\xi$ and $\mathcal Y^{\xi, \zeta}$. 
\end{rem}

The main argument used in this section to formally deduce the form of the limit averaged equation is the same that one would use in the classical case, when the fast process is not of McKean-Vlasov type. The only difference is truly in the type of frozen process that emerges from carrying out the formal argument, with consequence on the type of ergodic limit one needs to study and then exploit. 

Finally, before moving on, we would like to remark that, because of the discussion in the introduction around the conditional law $\varsigma_t^{x, \rho}$,  it is natural to wonder whether there is some relation between the coarse grained process $\bar X$ we obtain in this paper and the coarse graining produced via so called "projection method" or "effective dynamics using conditional expectation". The method, introduced in \cite{legoll2010effective}, approximates SDEs of the form \eqref{eq:multiscale system intro}  by the following effective equation for the slow variable:
$$
d \widetilde{\mathcal X}_t = {v}_P( \widetilde{\mathcal X}_t) dt + \sqrt2 \kappa dB_t \, \quad v_P(x) := \int_{R^m} v(x,y) \mu_{\infty}(y\vert x) dy \, 
$$
where the measure  $\mu_{\infty}(y\vert x)$ is the equilibrium conditional measure of system \eqref{eq:multiscale system intro}. Namely, let $\mu_{\infty}$ be the equilibrium measure of system \eqref{eq:multiscale system intro} (assuming such a measure exists). Then  $\mu_{\infty}(y\vert x) = \mu_{\infty} (x,y)/ \int_{\R^d} \mu_{\infty} (x,y) dx$ . In general this measure would depend on $\eps$ but this fact is usually disregarded and the measure is typically calculated by taking $\eps=1$. The difference between averaging and the projection method (as well as other related matters) has been discussed at length in \cite{duong2025coarse}. 
Here we just remark that one easy way to see that the two procedures don't coincide is to consider the case when the process $Y^{\eps}$ in \eqref{eq:def:McKean-Vlasov_full} is decoupled from $X^{\eps}$. In that case the conditional measure would depend on $x$ while the procedure we use would boil down to classic averaging, as we have already pointed out, and hence the measure we use to "project against" would just be the equilibrium measure of the fast process, which does not depend on $x$.



\section{Assumptions and main result}
\label{s:Assumptions_main_result}

In what follows, we work under the following assumptions. Let $\Wass_2$ denote the Wasserstein--2 distance on $\mathcal{P}_2(\mathbb{R}^d)$ or $\mathcal{P}_2(\mathbb{R}^m)$ (the underlying space will always be clear from context and for a definition of $\Wass_2$ see next section). 

\begin{assu}
\label{assu:gen} We make the following assumptions.
    \begin{enumerate}
        \item 
        There exist constants $K_1, K_2, K_3 > 0$ such that 
        \[
        |v(x_1, y_1) - v(x_2, y_2)|^2 + |b(x_1, y_1, \nu_1) - b(x_2, y_2, \nu_2)|^2 
        \leq 
        K_1 |x_1 - x_2|^2 
        + K_2 |y_1 - y_2|^2
        + K_3 \Wass_2^2(\nu_1, \nu_2)
        \]
        for any $x_1, x_2 \in \R^d$, $y_1, y_2 \in \R^m$, and $\nu_1, \nu_2 \in \mathcal{P}_2(\R^m)$.\footnote{
        Note that our proofs  require Assumption \ref{assu:gen} $(1)$ for $v$ only. The same assumption on $b$ is imposed to guarantee the existence and uniqueness of the solution to \eqref{eq:def:McKean-Vlasov_full}, for which we use a known result (see Proposition~\ref{prop:existence_uniqueness}).
        }
        \item The drift function $v$ is bounded, i.e. there exists $M_v > 0$ such that
        \[
            |v(x, y)| \leq M_v
        \]
        for any $x \in \R^d$, $y \in \R^m$.
        \item There exist constants $\beta_1, \beta_2, \beta_3 > 0$ such that 
        \begin{equation}
        \label{eq:assu_on_contractivity}
            2 \big \langle y_1 - y_2, b(x_1, y_1, \nu_1) - b(x_2, y_2, \nu_2) \big\rangle  
            \leq 
            - \beta_1 |y_2 - y_1|^2 
            + \beta_2 |x_2 - x_1|^2 
            + \beta_3 \Wass_2^2(\nu_1, \nu_2).          
        \end{equation}
        Assume additionally that $\beta_1 > \beta_3$.
    \end{enumerate}
\end{assu}

\begin{rem} Some straightforward comments on the above set of assumptions.
\begin{enumerate}
    \item The first assumption gives that there exists $K_4 > 0 $ such that 
    \be\label{eqn: linear growth v}
    |v(x, y)|^2 \leq K_4 (1 + |x|^2 + |y|^2)
    \ee
    \item From \eqref{eq:assu_on_contractivity} it follows that there exist constants
    $\beta_0^\prime, \beta_1^\prime,  > 0$
    such that
    \begin{equation}
    \label{eq:rem:beta_prime}
        2 \big\langle y, b(x, y, \nu)\big\rangle
        \leq
        \beta_0^\prime
        - \beta_1^\prime |y|^2
        + \beta_2 |x|^2
        + \beta_3 \int |u|^2 \nu(\dd{u}),
    \end{equation}
    with $\beta_1^\prime > \beta_3$. Indeed, from \eqref{eq:assu_on_contractivity}, we have
    \begin{equation*}
        2 \big\langle y, b(x, y, \nu)\big\rangle
        \leq
        2 \big\langle y, b(0, 0, \delta_0)\big\rangle
        - \beta_1 |y|^2
        + \beta_2 |x|^2
        + \beta_3 \int |u|^2 \nu(\dd{u}).
    \end{equation*}
    By Young's inequality, for any $a > 0$, we have
    \[
        2 \big\langle y, b(0, 0, \delta_0) \big\rangle
        \leq
        a |y|^2 + \frac{1}{a}|b(0, 0, \delta_0)|^2.
    \]
    Therefore, \eqref{eq:rem:beta_prime} follows by choosing $a>0$
    sufficiently small and setting
    \[
        \beta_0^\prime = \frac{|b(0,0,\delta_0)|^2}{a},
        \qquad
        \beta_1^\prime = \beta_1 - a,
    \]
    In particular, since $\beta_1 > \beta_3$, one can choose
    $a \in (0,\beta_1-\beta_3)$ so that
    $\beta_1^\prime > \beta_3$.
    \end{enumerate}
\end{rem}

First, let us make sure that the processes that we study are well defined. Consider the following proposition, the proof of which was obtained in \cite{XLLM21}. The authors of this paper used more general assumptions than ours, so Proposition~\ref{prop:existence_uniqueness} is a restatement of \cite[Theorem 3.1 and Lemma 5.1]{XLLM21} in our context and with our notation.
\begin{prop}[{\cite[Theorem 3.1 and Lemma 5.1]{XLLM21}}]
    \label{prop:existence_uniqueness}
    Under Assumption~\ref{assu:gen}, there exists a unique strong solution $(X^\eps_t, Y^\eps_t, \nu^\eps_t)_{t\geq 0}$ to the system \eqref{eq:def:McKean-Vlasov_full}, where $(X^\eps_t)_{t\geq 0}$ and $(Y^\eps_t)_{t\geq 0}$ are $\R^d$ and $\R^m-$ valued stochastic processes, respectively, and $\nu_t^\eps = \Law(Y_t^\eps) \in \mathcal{P}_2(\R^m)$ for any $t \geq 0$.

    Additionally, the second moments of $X^\eps_t$ and $Y^\eps_t$ are bounded uniformly in $\eps$ for finite time intervals, i.e. for any $T > 0$, there exists a constant $\widetilde{C} > 0$ such that
     \begin{equation*}
         \sup_{0 < \eps < 1}\sup_{0 \leq t \leq T} \E |X^\eps_t|^2 < \widetilde{C}, \qquad 
         \sup_{0 < \eps < 1}\sup_{0 \leq t \leq T} \E |Y^\eps_t|^2 < \widetilde{C}.
     \end{equation*}
\end{prop}

The next result addresses the well-posedness of the stationary problem \eqref{eqn: Pi intro}.  
 
Let us recall that 
for any given $\rho \in \mathcal P_2$, we seek a family $(\Pi^{x, \rho})_{x \in \R^d}$ which is a solution to equation  \eqref{eqn: Pi intro}.
It turns out that, under our assumptions, not only does a unique solution to \eqref{eqn: Pi intro} exists, but it also has a smooth density and enjoys a Lipschitz property with respect to $x$ and $\rho$, as stated in the following proposition.
\begin{prop}[Existence and uniqueness of $\Pi^{x, \rho}$]
    \label{prop:Pi_x_rho_existence}
    Under Assumption~\ref{assu:gen}, for any $\rho \in \mathcal{P}_2(\R^d)$, there exists a unique family of measures $(\Pi^{x, \rho})_{x \in \R^d}$ in $\mathcal{P}_2(\R^m)$ that solves equation \eqref{eqn: Pi intro}; such measures  have smooth densities for any $x \in \R^d$, $\rho \in \mathcal{P}_2(\R^d)$.
    
     Moreover, this family is Lipschitz continuous in both parameters: there exist constants $\widetilde{K}_1, \widetilde{K}_2 > 0$ such that
    \[
        \Wass_2^2(\Pi^{x_1, \rho_1}, \Pi^{x_2, \rho_2}) \leq \widetilde{K}_1 |x_1 - x_2|^2 + \widetilde{K}_2 \Wass_2^2(\rho_1, \rho_2).
    \]
\end{prop}

Another result that plays an important role in our proof is a result on the ergodic behavior of the time average of the process $\mathcal Y^{\xi}$ introduced in \eqref{eqn: our frozin proces intro}, which from now on we denote $\mathcal{Y}^{\xi, \zeta}$, as it will be useful to keep track of the initial datum in the notation:  
\begin{equation}
\label{eq:frozen_process_1_gen}
\begin{cases}
    \dd{\mathcal{Y}^{\xi, \zeta}_t} = b (\xi, \mathcal{Y}^{\xi, \zeta}_t, \boldsymbol{\nu}^{\xi, \zeta}_t) \dd{t} + \sqrt{2} \dd{\widetilde{W}_t}, \\
    \boldsymbol{\nu}^{\xi, \zeta}_t = \Law(\mathcal{Y}^{\xi, \zeta}_t) \\
    \mathcal{Y}^{\xi, \zeta}_0 = \zeta \text{ a.s.}
\end{cases}    
\end{equation}
The random variables  $\xi \in L^2(\Omega;\R^d),\zeta \in L^2(\Omega;\R^m)$ can in general depend on each other but we require them to be jointly independent of the Brownian motion $\widetilde W$.  

\begin{prop}[Ergodicity of $\mathcal{Y}^{\xi,\zeta}$]
    \label{prop:ergodicity_of_Y}
    For any $\rho \in \mathcal{P}_2(\R^d)$, let $(\Pi^{z,\rho})_{z \in \R^d}$ be the probability measures
    defined as solutions of \eqref{eqn: Pi intro}. 
    Let  $f:\R^d\times\R^m\to\R$ be any  measurable function and assume and that there exist
    constants $L_f,M_f\ge 0$ such that
    \[
    |f(z,y)-f(z,y')|\le L_f |y-y'|,
    \qquad 
    |f(z,y)|\le M_f,
    \qquad z\in\R^d,\ \ y,y'\in\R^m .
    \]
    Then, for any $\xi \in L^2(\Omega; \R^d)$ and for any $\zeta\in L^2(\Omega;\R^m)$ jointly independent of the Brownian
    motion $\widetilde W$ driving equation~\eqref{eq:frozen_process_1_gen}, we have
    \begin{equation}
    \label{eq:as_ergodicity_Y}
    \left|
    \frac1t\int_0^t f(\xi,\mathcal Y_s^{\xi,\zeta})\,\dd s
    -
    \int_{\R^m} f(\xi,y)\,\Pi^{\xi,\rho}(\dd y)
    \right|
    \xrightarrow[t\to\infty]{\Prob\text{--a.s.}} 0,
    \end{equation}
    where $\rho = \Law(\xi)$. Moreover, the convergence also holds in $L^2(\Omega;\R)$, uniformly over
    sets of initial random variables with uniformly bounded second moments. That is, for any fixed $C>0$, we let 
    \begin{equation}
    \label{eq:def:K_C}
        \mathcal K_C
        :=
        \left\{
        \zeta \in L^2(\Omega;\R^m):
        (\xi, \zeta) \text{ is independent of } \widetilde W
        \text{ and }
        \E|\zeta|^2 \leq C
        \right\}. 
        \footnote{
            In principle, the set $\mathcal K_C$ depends on the choice of $\xi$. Yet, we don't indicate this dependence in our notation, since it will not play a major role in the proof.
        }
    \end{equation}
    Then, for any $\xi \in L^2(\Omega; \R^d)$, 
    \begin{equation}
    \label{eq:L2_ergodicity_Y}
    \lim_{t\to\infty}
    \sup_{\zeta \in \mathcal{K}_C}
    \E\!\left[
    \left|
    \frac1t\int_0^t f(\xi,\mathcal Y_s^{\xi,\zeta})\,\dd s
    -
    \int_{\R^m} f(\xi,y)\,\Pi^{\xi,\rho}(\dd y)
    \right|^2
    \right]
    =0 ,
    \end{equation}
     and in the above the expectation is taken with respect to the underlying probability triplet $(\Omega, \mathcal F, \Prob )$  i.e. with respect to all sources of  randomness. 
\end{prop}


\begin{rem}\label{note: convergence}
The limit in Proposition~\ref{prop:ergodicity_of_Y} is random and, consistently  with our observations in Note~\ref{rem:Y_xi_zeta_parameters}, the almost sure convergence stated in the above proposition should be intended as follows:   for $\Prob$--a.e.\ $\omega$,
\[
    \left|
    \frac1t\int_0^t f\bigl(\xi(\omega),\mathcal Y_s^{\xi,\zeta}(\omega)\bigr)\,\dd s
    -
    \int_{\R^m} f\bigl(\xi(\omega),y\bigr)\,\Pi^{\xi(\omega),\rho}(\dd y)
    \right|
    \xrightarrow[t\to\infty]{}0 .
\]
Studying the above convergence boils down to understanding the behaviour of the conditional law $\varsigma^{x, \rho}$ solution of \eqref{eqn: evol conditional law intro}, and in particular the fact that it converges, as $t \rightarrow \infty$,   to $\Pi^{x, \rho}$. To see this, we can just observe that the expectation in \eqref{eq:L2_ergodicity_Y} can be rewritten, using the joint law  $\mathcal L(\xi, \mathcal Y^{\xi, \zeta}_t) = \boldsymbol{\mu}^{\xi, \zeta}_t = \rho \cdot \varsigma^{x,\rho}$, as 
\begin{align*}
    & \mathbb E \left(\frac{1}{t} \int_0^t f(\xi, \mathcal Y^{\xi, \zeta}_s)) \dd{s}
    - \int_{\R^m} f(\xi, y) \Pi^{\xi, \rho}(y) \dd{y} \right)  \\
    & = 
    \frac{1}{t} \int_0^t \dd{s}\left[\int_{\R^d} \rho(\dd{x}) \int_{\R^m}  \dd{y} f(x,y) \varsigma_s^{x, \rho}(y) \right] - \int_{\R^d} \rho(\dd{x}) \int_{\R^m} f(x, y) \Pi^{x, \rho}(y) \dd{y}.
\end{align*}
We don't prove explicitly the convergence of $\varsigma^{x, \rho}$ to $\Pi^{x,\rho}$ but, in view of the above observation, this is a consequence of Proposition \ref{prop:ergodicity_of_Y}.
\end{rem}

Now we are ready to present the main result. 
\begin{thm}[Main result]
\label{thm:main_thm}
    Let $X^\eps$ be as in \eqref{eq:def:McKean-Vlasov_full} and $\bar{X}$ be as in \eqref{eq:X_bar_rho_bar}. Let Assumption~\ref{assu:gen} hold. 
    Then the slow process $X^\eps$ converges to $\bar{X}$ in $L^2(\Omega; \R^d)$, uniformly on finite time intervals, i.e. for any $T > 0$, we have
    \begin{equation*}
        \lim_{\eps \to 0} \sup_{0\leq t \leq T} \E |X_t^\eps - \bar{X}_t|^2 = 0.
    \end{equation*}
\end{thm}

The proof of Theorem~\ref{thm:main_thm} is provided in Section~\ref{s:proof_main_result}.

\section{Proof of Proposition \ref{prop:Pi_x_rho_existence} and of Proposition \ref{prop:ergodicity_of_Y}}\label{s: proofs auxiliary results}
In this section we prove  Proposition \ref{prop:Pi_x_rho_existence} and of Proposition \ref{prop:ergodicity_of_Y}, starting from the latter, which is less standard. 
Throughout the paper, we use the following standard terminology. Given
$\mu,\nu\in\mathcal P_2(\R^k)$, a coupling of $\mu$ and $\nu$ is a probability
measure $\gamma$ on $\R^k\times\R^k$, such that $\int \gamma(\cdot, \dd{y}) = \mu$ and $\int \gamma(\dd{x}, \cdot) = \nu$, i.e. the first and second marginals are $\mu$ and
$\nu$, respectively. The quadratic Wasserstein distance between $\mu$ and
$\nu$ is defined by
\[
    \Wass_2^2(\mu,\nu)
    :=
    \inf_{\gamma}
    \int_{\R^k\times\R^k}
    |x-y|^2\,\gamma(\dd x,\dd y),
\]
where the infimum is taken over all couplings $\gamma$ of $\mu$ and $\nu$.
This infimum is always attained, and any coupling achieving it is called an
optimal coupling; see, for instance, \cite[Theorem~4.1]{villani2009optimal}.
Also, by definition of the $\Wass_2$ distance, it trivially follows that for any random variables $Z_1 \sim \mu$ and $Z_2 \sim \nu$, the following inequality holds
\begin{equation}
\label{eq:Wass_2_property}
    \Wass_2^2(\mu, \nu) \leq \E |Z_1 - Z_2|^2,
\end{equation}
which will be used multiple times in the proof.

\subsection{Proof of Proposition~\ref{prop:ergodicity_of_Y}}\label{subsec: proof of prop 3.5} We start by presenting the logic of the proof, the actual proof is detailed just afterwards. 

Recall the notation
\[
    \boldsymbol\mu_t^{\xi,\zeta}(\dd z,\dd y)
    =
    \Law\bigl(\xi,\mathcal Y_t^{\xi,\zeta}\bigr)(\dd z,\dd y),
\]
for the joint law of
\(\bigl(\xi,\mathcal Y_t^{\xi,\zeta}\bigr)\); then
\(\boldsymbol\mu^{\xi,\zeta}\) is a solution to the following system of
integro-differential equations:
\begin{equation}
\label{eq:def:mu_xi_zeta}
\begin{cases}
    \partial_t \boldsymbol\mu_t^{\xi,\zeta}
    =
    -\nabla_y\cdot
    \bigl(
        b(x,y,\boldsymbol\nu_t^{\xi,\zeta})
        \boldsymbol\mu_t^{\xi,\zeta}
    \bigr)
    +
    \Delta_y\boldsymbol\mu_t^{\xi,\zeta},
    \\[1mm]
    \boldsymbol\nu_t^{\xi,\zeta}(\dd y)
    =
    \displaystyle
    \int_{\R^d}
    \boldsymbol\mu_t^{\xi,\zeta}(\dd x,\dd y).
\end{cases}
\end{equation}
In particular,
\[
    \boldsymbol\nu_t^{\xi,\zeta}
    =
    \Law\bigl(\mathcal Y_t^{\xi,\zeta}\bigr).
\]
Note that the unknown \(\boldsymbol\mu_t^{\xi,\zeta}\) is a measure on
\(\R^d\times\R^m\), and \(x\) enters the first equation not as a parameter,
as in the case of the frozen process in classical averaging, but as a
variable. Indeed, by the second equation,
\(\boldsymbol\nu_t^{\xi,\zeta}\) is determined by
\(\boldsymbol\mu_t^{\xi,\zeta}\) and depends on the values of \(\boldsymbol\mu_t^{\xi,\zeta}\) for different
values of \(x\), since it is obtained by integrating out the first variable.
In other words, \eqref{eq:def:mu_xi_zeta} is an integro-differential system:
it is differential in the \(y\)-variable and involves integration with
respect to the \(x\)-variable. Set
\begin{equation}
\label{eq:def:mu_tilde_nu_tilde}
    \widetilde{\boldsymbol\mu}^{\rho}(\dd x,\dd y)
    :=
    \rho(\dd x)\Pi^{x,\rho}(\dd y),
    \qquad
    \widetilde{\boldsymbol\nu}^{\rho}(\dd y)
    :=
    \int_{\R^d}
    \rho(\dd x)\Pi^{x,\rho}(\dd y).
\end{equation}
By the definition of \(\Pi^{x,\rho}\), see
Equation~\eqref{eqn: Pi intro},
\(\widetilde{\boldsymbol\mu}^{\rho}\) is a stationary solution to
\eqref{eq:def:mu_xi_zeta}, and
\(\widetilde{\boldsymbol\nu}^{\rho}\) is its second marginal. As will be
shown later in the proof, we in fact have
\[
    \boldsymbol\mu_t^{\xi,\zeta}
    \longrightarrow
    \widetilde{\boldsymbol\mu}^{\rho}
    \quad\text{and}\quad
    \boldsymbol\nu_t^{\xi,\zeta}
    \longrightarrow
    \widetilde{\boldsymbol\nu}^{\rho}
    \quad\text{as }t\to\infty,
\]
which will be an essential part of the proof of the ergodic limit.

There are two main issues in studying the ergodic limit of
\(\mathcal Y^{\xi,\zeta}\). First, the drift depends on the random parameter
\(\xi\), and second, the process is nonlinear.

The first issue implies that the equation for
\(\boldsymbol\nu^{\xi,\zeta}\) is not closed (as can be seen by formally integrating over $x$ equation \eqref{eq:def:mu_xi_zeta}), unlike the equation for the
joint law \(\boldsymbol\mu^{\xi,\zeta}\). This suggests studying the process
\[
    t\mapsto
    \bigl(\xi,\mathcal Y_t^{\xi,\zeta}\bigr)
\]
instead of \(\mathcal Y^{\xi,\zeta}\) alone.

To circumvent the second issue, we introduce a process \(Z\) solving the
following SDE:
\begin{equation}
\label{eq:def:process_Z}
\begin{cases}
    \dd Z_t
    =
    b\bigl(
        \xi,
        Z_t,
        \widetilde{\boldsymbol\nu}^{\rho}
    \bigr)\dd t
    +
    \sqrt{2}\,\dd\widetilde W_t,
    \\[1mm]
    (\xi,Z_0)\text{ is independent of }\widetilde W
    \quad\text{and}\quad
    \Law(\xi,Z_0)
    =
    \widetilde{\boldsymbol\mu}^{\rho}.
\end{cases}
\end{equation}
Thus, the equation for \(Z\) is obtained from that of
\(\mathcal Y^{\xi,\zeta}\) by replacing the time-dependent measure argument
\(\boldsymbol\nu_t^{\xi,\zeta}\) with its expected long-time limit
\(\widetilde{\boldsymbol\nu}^{\rho}\), while retaining the random first
argument \(\xi\). The initial condition is chosen so that
\[
    \Law(\xi,Z_0)
    =
    \widetilde{\boldsymbol\mu}^{\rho}.
\]
Consequently, the joint process
\(t\mapsto(\xi,Z_t)\) is stationary. An equivalent condition would be
\begin{equation}
\label{eq:conditional_of_Z_0}
    \Law(Z_0\mid\xi=x)
    =
    \Pi^{x,\rho}
\end{equation}
for \(\rho\)-almost every \(x\).
We emphasize that in the above \(\boldsymbol\nu_t^{\xi,\zeta}\) is given, so the process $Z$ is linear. 
Freezing the nonlinear measure dependence in the drift of a
McKean--Vlasov process at its unique stationary measure is expected, under
suitable conditions on the drift, to produce a linear process with the same
long-time behaviour as the original nonlinear dynamics; see, for instance,
\cite{pavliotis2025linearization}. In our setting, however, the drift of
\(\mathcal Y^{\xi,\zeta}\) depends on the random variable \(\xi\), unlike
in the standard case of McKean--Vlasov processes, so this result cannot be used directly.

Another useful way to interpret the definition of \(Z\), which is also in
the spirit of our proof, is to compare the two joint processes
\[
    t\mapsto
    \bigl(\xi,\mathcal Y_t^{\xi,\zeta}\bigr)
    \qquad\text{and}\qquad
    t\mapsto(\xi,Z_t).
\]
From this perspective, the latter is precisely the linearisation of the
former obtained by replacing its evolving joint law $(\boldsymbol{\mu}^{\xi, \zeta}_t)_{t \geq 0}$ with the stationary distribution
\(\widetilde{\boldsymbol\mu}^{\rho}\). Owing to the particular form of the
drift \(b\), only the marginal
\(\widetilde{\boldsymbol\nu}^{\rho}\) enters the resulting SDE for \(Z\).

Since the process \(t\mapsto(\xi,Z_t)\) is stationary, we can apply
Birkhoff's ergodic theorem to it. Provided that we prove that \(Z\) and
\(\mathcal Y^{\xi,\zeta}\) have the same long-time behaviour, we can then
transfer the resulting ergodic limit to the latter.

Birkhoff's ergodic theorem applied to
\(t\mapsto(\xi,Z_t)\) gives
\[
    \frac{1}{t}\int_0^t f(\xi,Z_s)\,\dd s
    \xrightarrow[t\to\infty]{}
    \E\bigl[f(\xi,Z_0)\mid\mathcal I\bigr]
    \quad \text{a.s. and in } L^1,
\]
where \(\mathcal I\) is the invariant $\sigma$-algebra, that is, the
$\sigma$-algebra consisting of events that are invariant under time shifts.
The precise definition will be given in the proof below.\footnote{
    In the proof, we denote by $\mathcal{I}$ the invariant $\sigma$-algebra on the space of trajectories, rather than on the probability space $(\Omega, \Prob)$, which doesn't change the argument
    }
Since \(f\) is bounded,
the convergence also holds in \(L^2\).

For an ergodic stationary diffusion, the invariant $\sigma$-algebra is
trivial. However, this is not the case for the process
\(t\mapsto(\xi,Z_t)\), since its first component does not change in time. It turns out that the invariant sigma algebra \(\mathcal I\) coincides with the sigma algebra generated by $\xi$, $\mathcal I=\sigma(\xi)$. We will prove this rigorously below, although, upon reflection,  this fact is quite intuitive. One way of understanding this is by observing that, since the first component of $t \mapsto (\xi, Z_t)$ is fixed in time, and, 
conditionally on $\xi = x$, the second component is ergodic, it is necessary and sufficient
to know the outcome of $\xi$ in order to determine the value of the time limit 
\[
    \frac{1}{t} \int_0^t f(\xi, Z_s)\dd{s}.
\]
When we say that,  conditionally on $\xi = x \in \R^d$, the process $Z$ is ergodic, we mean that, conditionally on $\xi = x \in \R^d$, such a process behaves like 
the linear stationary process $Z^x$ defined as: 
\begin{equation}
    \label{eq:def:Z_x}
    \begin{cases}
        \dd{Z_t^x} = b(x, Z_t^x, \widetilde{\boldsymbol{\nu}}^\rho)\dd{t} + \sqrt{2}\dd{\widetilde{W}_t}, \\
        Z_0^x \sim \Pi^{x, \rho},
    \end{cases}
\end{equation}
where the initial condition comes from \eqref{eq:conditional_of_Z_0}.
Under our assumptions, this process is, in fact, ergodic (in the sense that the time average
is equal to the ensemble average with respect to the unique stationary 
distribution $\Pi^{x, \rho}$).

By \eqref{eq:conditional_of_Z_0}, the identification $\mathcal I = \sigma(\xi)$ leads to:
\[
    \E\bigl[f(\xi,Z_0)\mid\mathcal I\bigr]
    =
    \E\bigl[f(\xi,Z_0)\mid\xi\bigr]
    =
    \int_{\R^m}
    f(\xi,y)\Pi^{\xi,\rho}(\dd y) \quad \text{a.s.},
\]
which proves Proposition~\ref{prop:ergodicity_of_Y}, provided that we could show that $\mathcal{Y}^{\xi, \zeta}$ has the same long time limit as $Z$.

This motivates the following proof.

\begin{proof}[Proof of Proposition \ref{prop:ergodicity_of_Y}]
The proof is in three steps. \\
\textbf{Step 1.}
Prove the following inequality:
\begin{equation}
\label{eq:Y-Z_as_comparison}
    \bigl|
        \mathcal Y_t^{\xi,\zeta}-Z_t
    \bigr|^2
    \leq
    \left[
        |\zeta-Z_0|^2
        +
        \E|\zeta-Z_0|^2
    \right]
    \e^{-(\beta_1-\beta_3)t}
    \quad\text{a.s.}
\end{equation}
In particular, this implies that
\(\mathcal Y_t^{\xi,\zeta}-Z_t\) converges to zero almost surely at an
exponential rate. Moreover, it follows that
\begin{equation}
\label{eq:Y-Z_L2_comparison}
    \sup_{\zeta\in\mathcal K_C}
    \E\bigl|
        \mathcal Y_t^{\xi,\zeta}-Z_t
    \bigr|^2
    \leq
    2
    \sup_{\zeta\in\mathcal K_C}
    \E|\zeta-Z_0|^2
    \e^{-(\beta_1-\beta_3)t},
\end{equation}
where the constant in front of the exponential is finite by the definition
of \(\mathcal K_C\) and the finite second moment of
\(\widetilde{\boldsymbol\mu}^{\rho}\).

\textbf{Step 2.}
Use Birkhoff's ergodic theorem to prove the ergodic result for \(Z\) in the
form
\begin{equation}
\label{eq:ergodic_result_for_Z}
    \frac{1}{t}\int_0^t f(\xi,Z_s)\,\dd s
    \xrightarrow[t\to\infty]{}
    \int_{\R^m}
    f(\xi,y)\Pi^{\xi,\rho}(\dd y)
    \quad\text{a.s. and in }L^2.
\end{equation}

\textbf{Step 3.}
Combine the first two steps to prove the proposition. The \(L^2\)-convergence
uniformly over \(\zeta\in\mathcal K_C\) will follow from
\eqref{eq:Y-Z_L2_comparison}, together with the fact that the definition of the
stationary process \(t\mapsto(\xi,Z_t)\) does not involve  \(\zeta\).

Now, let us provide complete proofs of each of the steps. 

{\bf Proof of Step 1.}

Since $\mathcal Y^{\xi,\zeta}$ and $Z$ are driven by the same Brownian
motion, their stochastic terms cancel. Hence, almost surely,
\[
    \frac{\dd}{\dd t}|\mathcal{Y}_t^{\xi,\zeta} - Z_t|^2
    =
    2\Bigl\langle
        \mathcal{Y}_t^{\xi,\zeta} - Z_t,\,
        b\bigl(
            \xi,\mathcal Y_t^{\xi,\zeta},
            \boldsymbol\nu_t^{\xi,\zeta}
        \bigr)
        -
        b\bigl(
            \xi,Z_t,\widetilde{\boldsymbol{\nu}}^\rho
        \bigr)
    \Bigr\rangle
    \leq
    -\beta_1|\mathcal{Y}_t^{\xi,\zeta} - Z_t|^2
    +
    \beta_3
    \Wass_2^2\bigl(
        \boldsymbol\nu_t^{\xi,\zeta},
        \widetilde{\boldsymbol{\nu}}^\rho
    \bigr).
\]
Taking expectations, by \eqref{eq:Wass_2_property}, we obtain
\[
\begin{aligned}
    \frac{\dd}{\dd t}\E|\mathcal{Y}_t^{\xi,\zeta} - Z_t|^2
    &\leq
    -\beta_1\E|\mathcal{Y}_t^{\xi,\zeta} - Z_t|^2
    +
    \beta_3
    \Wass_2^2\bigl(
        \boldsymbol\nu_t^{\xi,\zeta},
        \widetilde{\boldsymbol{\nu}}^\rho
    \bigr)
    \leq
    -(\beta_1-\beta_3)\E|\mathcal{Y}_t^{\xi,\zeta} - Z_t|^2.
\end{aligned}
\]
Therefore, by Gronwall's inequality,
\begin{equation}
\label{eq:expected_YZ_comparison}
    \E|\mathcal{Y}_t^{\xi,\zeta} - Z_t|^2
    \leq
    \E|\zeta-Z_0|^2
    \e^{-(\beta_1-\beta_3)t}.
\end{equation}
Consequently, by \eqref{eq:Wass_2_property},
\[
    \Wass_2^2\bigl(
        \boldsymbol\nu_t^{\xi,\zeta},
        \widetilde{\boldsymbol{\nu}}^\rho
    \bigr)
    \leq
    \E|\zeta-Z_0|^2\e^{-(\beta_1-\beta_3)t}.
\]

Returning to the pathwise differential inequality and using the preceding
estimate, we obtain
\[
    \frac{\dd}{\dd t}|\mathcal{Y}_t^{\xi,\zeta} - Z_t|^2
    \leq
    -\beta_1|\mathcal{Y}_t^{\xi,\zeta} - Z_t|^2
    +
    \beta_3 \E|\zeta-Z_0|^2
    \e^{-(\beta_1-\beta_3)t}.
\]
An application of the variation-of-constants formula gives
\[
\begin{aligned}
    |\mathcal{Y}_t^{\xi,\zeta} - Z_t|^2
    &\leq
    \e^{-\beta_1t}|\zeta-Z_0|^2
    +
    \beta_3\E|\zeta-Z_0|^2
    \int_0^t
    \e^{-\beta_1(t-s)}
    \e^{-(\beta_1-\beta_3)s}\dd s
    \\
    &=
    \e^{-\beta_1t}|\zeta-Z_0|^2
    +
    \E|\zeta-Z_0|^2
    \left(
        \e^{-(\beta_1-\beta_3)t}
        -
        \e^{-\beta_1t}
    \right)
    \\
    &\leq
    \bigl(
        |\zeta-Z_0|^2+\E|\zeta-Z_0|^2
    \bigr)
    \e^{-(\beta_1-\beta_3)t},
    \qquad \Prob\text{--a.s.}
\end{aligned}
\]
This proves \eqref{eq:Y-Z_as_comparison} and \eqref{eq:Y-Z_L2_comparison} after taking expectations on both sides.

{\bf Proof of Step 2.}
We now prove the ergodic result for the stationary process
$(\xi, (Z_t)_{t\geq 0})$. Most of what is in this step is substantially checking technical requirements in order to apply  Birkhoff's
ergodic theorem. 

Since $(\xi, (Z)_{t \geq 0})$ is stationary, it admits a two-sided
stationary extension $(\xi, (Z_t)_{t\in\R})$, that we will simply denote as $(\xi, Z)$. Consider the space
\[
    \mathcal{E}:= \R^d \times C(\R;\R^m),
\]
equipped with the product of the Euclidean topology on $\R^d$ and the topology of uniform convergence on compact subsets of $\R$ on $C(\R;\R^m)$.
Define the canonical embedding of the process into the space of
trajectories by
\[
    \Phi:\Omega\longrightarrow\mathcal{E},
    \qquad
    \Phi(\omega)
    :=
    (\xi(\omega),Z(\omega)).
\]
Since the trajectories of the process are $\Prob$-a.s.\ continuous,
after modifying the process on a $\Prob$-null set if necessary, the
mapping $\Phi$ is measurable.

The probability law $\mathbf P$ induced by the process
$(\xi,Z)$ on the path space $\mathcal E$ is defined by
\[
    \mathbf P:=\Prob\circ\Phi^{-1}.
\]
Equivalently,
\[
    \mathbf P(B)
    =
    \Prob(\Phi\in B),
    \qquad
    B \in \mathcal B(\mathcal E),
\]
where $\mathcal B(\mathcal E)$ denotes the Borel $\sigma$-algebra of
$\mathcal E$. We therefore constructed the probability space
$(\mathcal E,\mathcal B(\mathcal E),\mathbf P)$. Note that any set $B \in \mathcal{B}(\mathcal{E})$ can be decomposed as follows
\[
    B = \bigcup_{a \in A} \{a\} \times D_a, 
\]
for some $A \in \mathcal{B}(\R^d)$ and $D_a \in \mathcal{B}(C(\R; \R^m))$.
Furthermore,  $\rho = \Law(\xi)$ is precisely the marginal law of the first component on $\mathcal{E}$:
\[
    \rho (\cdot) = \int_{z \in C(\R; \R^m)} \mathbf P(\cdot, \dd{z}). 
\]
By the disintegration theorem, there exists a family of regular conditional distributions $(Q_x)_x$ on $C(\R; \R^m)$, for $\rho$-a.e. $x \in \R^d$, such that
\begin{equation}
\label{eq:desintegration_theorem}
     \mathbf P \Big(\bigcup_{a \in A} \{a\} \times D_a \Big) 
    = \int_{x \in A} \rho(\dd{x}) \int_{z \in D_a} Q_x(\dd{z}),
\end{equation}
for any $A \in \mathcal{B}(\R^d)$ and $D_a \in \mathcal{B}(C(\R; \R^d))$. 
Equivalently, $Q_x$ can be defined as a pull-back of the conditional probability 
$\Prob_x := \Prob (\cdot | \xi = x)$, i.e.:
\[
    Q_x = \Prob_x \circ \Phi^{-1}.
\]
Under the probability measure $\Prob_x$, the process $Z$ is $\Prob_x$-a.s. equal to the process $Z^x$ defined in \eqref{eq:def:Z_x}
Note that the process $Z^x$ is also stationary (both under $\Prob$ and $\Prob_x$). This directly follows from the definition of $\Pi^{x, \rho}$, $\widetilde{\boldsymbol{\nu}}^\rho$ (Equations~\eqref{eqn: Pi intro} and \eqref{eq:def:mu_tilde_nu_tilde}) and the Fokker-Planck equation for \eqref{eq:def:Z_x}. In particular, it means that its invariant $\sigma$-algebra
\[
    \mathcal{I}_x := \{D \in \mathcal{B}(C(\R; \R^m)): Q_x \big((\theta^{-1}_t D) \triangle D \big) = 0 \quad
    \text{for any } t \in \R\}
\]
is trivial, where $\triangle$ denotes the symmetric difference of two sets and $\theta_t$ is the time shift operator defined as follows.
For $t\in\R$, let
\[
    \theta_t:\mathcal E\longrightarrow\mathcal E,
    \qquad
    \theta_t (a, \varphi) := (a, \varphi_{t + \cdot}),
    \quad (a, \varphi) \in  \mathcal{E}.
\]
The stationarity of $(\xi,Z)$ implies that $\mathbf P$ is invariant
under every $\theta_t$, that is,
\[
    \mathbf P\circ\theta_t^{-1}=\mathbf P,
    \qquad t\in\R.
\]
With abuse of notation, we will also denote by $\theta_t$ the time-shift operator on $C(\R; \R^m)$.
Let
\[
    \mathcal I
    :=
    \left\{
        B\in\mathcal B(\mathcal E):
        \quad \mathbf{P} \big((\theta_t^{-1}B) \triangle B \big) = 0 \quad \text{for every }t\in\R
    \right\}
\]
be the invariant $\sigma$-algebra. Note that, by \eqref{eq:desintegration_theorem}, for any set $B \in \mathcal{I}$, we get:
\[
    0 = \mathbf{P} \big((\theta_t^{-1}B) \triangle B \big)
    =  \int_{x \in A} \rho(\dd{x}) \int_{z \in (\theta_t^{-1}D_x )\triangle D_x} Q_x(\dd{z}),
\]
which means that, for $\rho$-a.e. $x$, we have
\[
    \int_{z \in (\theta_t^{-1}D_x) \triangle D_x} Q_x(\dd{z}) = 0
\]
and, therefore, $D_x \in \mathcal{I}_x$. As we said earlier, $\mathcal{I}_x$ is trivial. 
Therefore, we get that, for any set $B \in \mathcal{I}$ of positive $\mathbf{P}$-measure, 
its second component is the full set $C(\R; \R^m)$, while the first component can be any $A \in \mathcal{B}(\R^d)$. In other words,
\[
    \Phi^{-1}\mathcal{I} = \sigma(\xi) \quad \Prob-\text{a.s.},
\]
which follows from the definition of the latter.

Finally, by a continuous-time ergodic theorem (see e.g. \cite[Corollary 25.9]{K21}), we get:
\begin{equation}
\label{eq:continuous_Birkhoff}
    \frac{1}{t} \int_0^t \psi (\theta_s(a, \varphi))\,\dd s
    \xrightarrow[t\to\infty]{}
    \mathbf E_{\mathbf P}
    [\psi \mid \mathcal I],
    \qquad
    \mathbf P\text{-a.s. and in }L^1(\mathcal{E}, \mathbf P),
\end{equation}
for any $\psi \in L^1(\mathcal{E}, \mathbf P)$.

Now, since $f(\xi,Z_0)\in L^1(\Omega, \Prob)$, we can define
\[
    \psi(a, \varphi) := f(a, \varphi(0)),
    \qquad (a, \varphi) \in \mathcal E.
\]
Then $\psi \in L^1(\mathcal E,\mathbf P)$ and 
\[
    \psi(\theta_s\Phi)
    =
    f(\xi,Z_s).
\]
Therefore, composing \eqref{eq:continuous_Birkhoff} with $\Phi$ gives 
\begin{equation*}
    \frac1t\int_0^t f(\xi,Z_s)\,\dd s
    \xrightarrow[t\to\infty]{}
    \mathbf E_{\mathbf P}
    [\psi \mid \mathcal I]\circ\Phi
    = \E
    [\psi \circ \Phi \mid \Phi^{-1}\mathcal I]
    = \E [f(\xi, Z_0) \mid \xi],
\end{equation*}
$\Prob\text{-a.s. and in }L^1(\Prob)$,
where the first equality follows from the definition of measure $\mathbf{P}$ and uniqueness of the conditional expectation.\footnote{
    To be more precise, if we denote $g = \mathbf{E}_{\mathbf{P}}[\psi | \mathcal{I}]$, then $g \circ \Phi$ is $\Phi^{-1} \mathcal{I}$ measurable and for any $B \in \mathcal{I}$, we have
    \begin{equation*}
        \begin{aligned}
            \int_{\Phi^{-1}(B)} g \circ \Phi \dd{\Prob} 
            = \int_{B} g \dd{\mathbf{P}}
            = \int_B \psi \dd{\mathbf{P}}
            = \int_{\Phi^{-1}(B)} \psi \circ \Phi^{-1} \dd{\Prob}
            = \int_{\Phi^{-1}(B)} \E [\psi \circ \Phi^{-1}| \Phi^{-1} \mathcal{I}] \dd{\Prob},
        \end{aligned}
    \end{equation*}
    where we used the properties of conditional expectation and change of variables. 
    Since the conditional expectation is unique, the identification 
    $g \circ \Phi = \E [\psi \circ \Phi^{-1}| \Phi^{-1} \mathcal{I}]$ follows $\Prob$-a.s.
    }
This proves the result \eqref{eq:ergodic_result_for_Z}, if we additionally use that the conditional 
law $\Law(Z_t | \xi = x) = \Pi^{x, \rho}$ for $\rho$-a.e. $x \in \R^d$.

Since $f$ is bounded,
\[
    \left|
        \frac1t\int_0^t f(\xi,Z_s)\dd s
        -
        \int_{\R^m}
        f(\xi,y)\Pi^{\xi,\rho}(\dd y)
    \right|^2
    \leq
    4M_f^2.
\]
Therefore, the almost-sure convergence in
\eqref{eq:ergodic_result_for_Z} and the dominated convergence theorem give
\begin{equation}
\label{eq:ergodic_Z_L2}
    \E\left[
        \left|
            \frac1t\int_0^t f(\xi,Z_s)\dd s
            -
            \int_{\R^m}
            f(\xi,y)\Pi^{\xi,\rho}(\dd y)
        \right|^2
    \right]
    \xrightarrow[t\to\infty]{}0.
\end{equation}
This proves \eqref{eq:ergodic_result_for_Z} in $L^2$.

\noindent
{\bf Proof of Step 3.}
We first prove the almost-sure convergence. By the Lipschitz continuity
of $f$ in its second variable,
\[
    \left|
        \frac1t\int_0^t
        f\bigl(
            \xi,\mathcal Y_s^{\xi,\zeta}
        \bigr)\dd s
        -
        \frac1t\int_0^t
        f(\xi,Z_s)\dd s
    \right|
    \leq
    \frac{L_f}{t}
    \int_0^t
    |
        \mathcal Y_s^{\xi,\zeta}-Z_s
    |\dd s.
\]
Therefore, using \eqref{eq:Y-Z_as_comparison} we obtain
\[
    \frac{L_f}{t}
    \int_0^t
    |
        \mathcal Y_s^{\xi,\zeta}-Z_s
    |\dd s
    \leq
    \frac{L_f}{t}
    \left(
        |\zeta-Z_0|^2+\E|\zeta-Z_0|^2
    \right)^{1/2}
    \int_0^\infty
        \e^{-(\beta_1-\beta_3)s/2}\dd s
    \xrightarrow[t\to\infty]{}0
\]
almost surely. Combining this with \eqref{eq:ergodic_result_for_Z} yields
\eqref{eq:as_ergodicity_Y}.

It remains to prove the $L^2$-convergence uniformly over
$\zeta\in\mathcal K_C$. For each $\zeta\in\mathcal K_C$, write
\[
\begin{aligned}
    &
    \frac1t\int_0^t
    f\bigl(
        \xi,\mathcal Y_s^{\xi,\zeta}
    \bigr)\dd s
    -
    \int_{\R^m}
    f(\xi,y)\Pi^{\xi,\rho}(\dd y)
    \\
    &\quad=
    \frac1t\int_0^t
    \left[
        f\bigl(
            \xi,\mathcal Y_s^{\xi,\zeta}
        \bigr)
        -
        f(\xi,Z_s)
    \right]\dd s
    +
    \frac1t\int_0^t
    \bigl[
        f(\xi,Z_s)-
        \int_{\R^m}
        f(\xi,y)\Pi^{\xi,\rho}(\dd y)
    \bigr]\dd s.
\end{aligned}
\]
By the Lipschitz continuity of $f$, and
\eqref{eq:Y-Z_L2_comparison},
\[
\begin{aligned}
    \sup_{\zeta\in\mathcal K_C}
    \left(
        \E\left|
            \frac1t\int_0^t
            \left[
                f\bigl(
                    \xi,\mathcal Y_s^{\xi,\zeta}
                \bigr)
                -
                f(\xi,Z_s)
            \right]\dd s
        \right|^2
    \right)^{1/2}
    &\leq
    \frac{L_f}{t}
    \int_0^t
    \sup_{\zeta\in\mathcal K_C}
    \left(
        \E|
            \mathcal Y_s^{\xi,\zeta}-Z_s
        |^2
    \right)^{1/2}\dd s
    \\
    &\leq
    \frac{L_fC^{1/2}}{t}
    \int_0^t
        \e^{-(\beta_1-\beta_3)s/2}\dd s
    \\
    &\leq
    \frac{
        2L_fC^{1/2}
    }{
        (\beta_1-\beta_3)t
    }
    \xrightarrow[t\to\infty]{}0.
\end{aligned}
\]

Combining this estimate with $L^2$ convergence in \eqref{eq:ergodic_result_for_Z} yields \eqref{eq:L2_ergodicity_Y}, which completes the proof.

\end{proof}

\subsection{Proof of Proposition~\ref{prop:Pi_x_rho_existence}}
Let us observe that  \eqref{eqn: Pi intro} can be equivalently written as the following system of equations: 
\begin{equation}
\label{eq:system_for_Pi_rho}
    \begin{cases}
        -\nabla_y\cdot
        \left[b(x,y,\nu)\Pi^{x,\rho}(y)\right]
        + \Delta_{y} \Pi^{x,\rho}(y) = 0, \\
        \nu(y) = \int \Pi^{z,\rho}(y)\,\rho(\dd z) \,.
    \end{cases}
\end{equation}

The proof will be given in 4 steps.

{\bf Step 1}: We start by solving the first of the equations in \eqref{eq:system_for_Pi_rho}, for given $x$ and $\nu$. That is,  we show that for any
$x\in\R^d$ and $\nu\in\mathcal P_2(\R^m)$ given, there exists a
unique measure $\pi^{x,\nu} \in \mathcal{P}_2(\R^m)$ solution to the following equation 
\begin{equation}
\label{eq:assu:unique_invariant_measure}
    -\nabla_y \cdot\Bigl[b(x,y,\nu)\pi^{x,\nu}(y)\Bigr]
    +\Delta_{y}\pi^{x,\nu}(y)
    =0.
\end{equation}

Moreover, since the above is a linear elliptic Fokker-Planck equation, the solution is a classical solution, which is the density of $\pi^{x, \nu}$ (and we keep denoting such a density by $\pi^{x, \nu}$ with abuse of notation).

{\bf Step 2}: Show that $\pi^{x,\nu}$ is Lipschitz continuous in $\nu$, i.e.
\[
    \Wass_2^2(\pi^{x, \nu}, \pi^{x, \nu^\prime}) 
    \leq C \Wass_2^2(\nu, \nu^\prime),    
\]
for some constant $C > 0$ (that can be explicitly determined) and for any $x \in \mathcal{P}_2(\R^d)$, $\nu, \nu^\prime \in \mathcal{P}_2(\R^m)$. 

{\bf Step 3}: As a result of Step 1, for each $\nu \in \mathcal P_2(\R^m)$ we have a one parameter family of measures indexed by $x \in \R^d$,  i.e. $\{\pi^{x, \nu}\}_{x}$. Here we show that, for any given $\rho \in \mathcal P_2(\R^d)$, there exists a unique measure $\nu^\rho \in \mathcal{P}_2(\R^m)$, such that $\nu^{\rho}$ (and the corresponding family $\{\pi^{x, \nu^{\rho}}\}_{x}$) satisfy  the following equation:
\[
    \nu^\rho(\dd y)
    =
    \int_{\R^d}\pi^{z,\nu^\rho}(\dd y)\rho(\dd z).
\]

{\bf Step 4}:
If we define
\[
    \Pi^{x,\rho} := \pi^{x,\nu^\rho} \in \mathcal{P}_2(\R^m),
\]
then $\Pi^{x,\rho}$ is indeed the unique solution of
\eqref{eq:system_for_Pi_rho},  by construction.
We finally show that $\Pi^{x,\rho}$ is Lipschitz continuous in both $x$ and $\rho$, i.e.
\[
    \Wass_2^2(\Pi^{x_1, \rho_1}, \Pi^{x_2, \rho_2}) 
    \leq 
    C_1 |x_1 - x_2|^2 + C_2 \Wass_2^2(\rho_1, \rho_2),
\]
for some constants $C_1, C_2 > 0$ (that will be specified) and for any $x_1, x_2 \in \R^d$, $\rho_1, \rho_2 \in \mathcal{P}_2(\R^d)$.

Now, let us provide complete proofs of each of the steps.

{\bf Proof of Step 1}. For fixed $x\in\mathbb{R}^d$ and
$\nu\in\mathcal{P}_2(\mathbb{R}^m)$, consider the diffusion on $\R^m$
\begin{equation}
    \label{eq:aux:Y_nu}
    \dd Y_t=b(x,Y_t,\nu)\,\dd t+\sqrt{2}\,\dd W_t,    
\end{equation}
whose infinitesimal generator is
\[
    L^{x,\nu}\varphi(y)
    =
    \langle b(x,y,\nu), \nabla\varphi(y) \rangle + \Delta\varphi(y).
\]
Note that $Y$ is a linear process whose stationary solution, if it exists, is precisely the solution to \eqref{eq:assu:unique_invariant_measure}. 
By Assumption~\ref{assu:gen}, the Lyapunov function
\[
    V(y):=1+|y|^2
\]
satisfies
\[
    (L^{x,\nu}V)(y)
    =
    2 \langle y, b(x,y,\nu) \rangle + 2m
    \leq
    -aV(y)+b_{x,\nu}
\]
for some $a>0$ and $b_{x,\nu}<\infty$. Consequently, upon setting
\[
    C_{x,\nu}
    :=
    \left\{
        y\in\mathbb{R}^m:
        V(y)\leq \frac{2b_{x,\nu}}{a}
    \right\},
\]
we obtain
\[
    (L^{x,\nu}V)(y)
    \leq
    -\frac{a}{2} V(y)
    +
    b_{x,\nu}\mathbf{1}_{C_{x,\nu}}(y).
\]
Thus, the Foster--Lyapunov condition \textnormal{(CD2)} of
\cite[Section~4.1]{meyn1993stability} is satisfied. Since the diffusion
Feller and $C_{x,\nu}$ is compact,
\cite[Theorem~4.5]{meyn1993stability} yields existance of an invariant
probability measure $\pi^{x,\nu}$ that satisfies
\[
    \int_{\mathbb{R}^m}
    V(y) \pi^{x,\nu}(\dd y)
    <\infty.
\]
In particular, it means that $\pi^{x,\nu}\in\mathcal{P}_2(\mathbb{R}^m)$. Moreover, the diffusion coefficient is constant and non-degenerate,
hence the invariant probability measure $\pi^{x,\nu}$ is unique and has a density.

{\bf Proof of Step 2}. We next establish the Lipschitz continuity of $\pi^{x,\nu}$ with respect to the parameter $\nu$. Let $Y$ be the solution of \eqref{eq:aux:Y_nu}, and let $Y^\prime$ be the solution of the same equation with $\nu$ replaced by $\nu^\prime\in\mathcal{P}_2(\mathbb{R}^m)$. We couple the two processes synchronously (using the same Brownian motion) and setting $Y_0=Y_0^\prime$ almost surely. 
Note that $\pi^{x, \nu}$ and $\pi^{x, \nu^\prime}$ are the unique stationary distributions of $Y$ and $Y^\prime$ respectively. Moreover, regardless of the starting law, the distributions $\Law(Y_t)$ and $\Law(Y^\prime_t)$ converge to their respective stationary measures in $\mathcal{P}_2(\R^m)$. For these two processes, by Assumption~\ref{assu:gen}, we have:
\begin{equation*}
\begin{aligned}
    \E|Y_t - Y_t^\prime|^2 &= 2 \int_0^t \E \big\langle Y_t - Y_t^\prime, b(x, Y_t, \nu) - b(x, Y_t^\prime, \nu^\prime) \big\rangle \dd{s} \leq - \beta_1 \int_0^t \E |Y_s - Y_s^\prime|^2 \dd{s} + \beta_3 \Wass_2^2(\nu, \nu^\prime).
\end{aligned}
\end{equation*}
From this inequality, it follows that
\begin{equation*}
    \sup_{t \geq 0} \E|Y_t - Y_t^\prime|^2 \leq  \frac{\beta_3}{\beta_1}\Wass_2^2(\nu, \nu^\prime),
\end{equation*}
which finally gives us: 
\begin{equation}
\label{eq:aux:Wass_pi_nu_nu_prime}
    \Wass_2^2(\pi^{x, \nu}, \pi^{x, \nu^\prime}) \leq \sup_{t \geq 0} \Wass_2^2(\Law(Y_t), \Law(Y_t^\prime)) \leq \sup_{t \geq 0} \E|Y_t - Y_t^\prime|^2 \leq \frac{\beta_3}{\beta_1}\Wass_2^2(\nu, \nu^\prime).
\end{equation}

{\bf Proof of Step 3}. The second equation of \eqref{eq:system_for_Pi_rho} can be understood as follows. For a $\rho \in \mathcal{P}_2(\R^d)$, define $\nu^\rho \in \mathcal{P}_2(\R^m)$ to be the solution to the following equation:
\begin{equation}
\label{eq:assu:unique_stationary_point}
    \nu^\rho(y) = \int \rho(\dd{z}) \pi^{z, \nu^\rho}(y),
\end{equation}
where $\pi^{z, \nu}$ is the unique solution of \eqref{eq:assu:unique_invariant_measure} for the respective $z \in \R^d$, $\nu \in \mathcal{P}_2(\R^m)$.
In other words, $\nu^\rho$ is the fixed point of the mapping $\Phi: \nu \mapsto \int \rho(\dd{x}) \pi^{x, \nu}(\cdot) $ on $\mathcal{P}_2(\R^m)$. In order to show the existence and uniqueness of such a fixed  point, we will show the contractivity of the mapping $\Phi$.

Take the optimal (in the sense of $\Wass_2$) coupling between $\pi^{x, \nu}$ and $\pi^{x, \nu^\prime}$, i.e. the measure $\gamma^x$ such that the infimum in the definition of the $\Wass_2$ distance is attained. Then, the following measure on $\R^m \times \R^m$:
\begin{equation*}
    \Gamma(A, B) := \int \rho(\dd{x}) \gamma^x(A, B),
\end{equation*}
where $A, B$ are some measurable subsets of $\R^m$, is a coupling between $\Phi(\nu)$ and $\Phi(\nu^\prime)$. Indeed, its marginals are
\begin{equation*}
    \int \Gamma(\dd{y}, \cdot) = \iint \gamma^x(\dd{y}, \cdot) \rho(\dd{x}) = \int \rho(\dd{x}) \pi^{x, \nu^\prime}(\cdot)  = \Phi(\nu^\prime), 
\end{equation*}
and, similarly, 
\begin{equation*}
    \int \Gamma(\cdot, \dd{y^\prime}) = \Phi(\nu).
\end{equation*}
Therefore, the Wasserstein distance between the two is bounded by:
\begin{equation*}
\begin{aligned}
    \Wass_2^2(\Phi(\nu), \Phi(\nu^\prime)) &\leq \int |y - y^\prime|^2 \Gamma(\dd{y}, \dd{y^\prime}) = \iint |y - y^\prime|^2 \gamma^x(\dd{y}, \dd{y^\prime}) \rho(\dd{x}) \\
    &= \int \Wass_2^2 (\pi^{x, \nu}, \pi^{x, \nu^\prime}) \rho(\dd{x}) \leq \frac{\beta_3}{\beta_1} \Wass_2^2 (\nu, \nu^\prime),    
\end{aligned}
\end{equation*}
where the last inequality holds due to \eqref{eq:aux:Wass_pi_nu_nu_prime}. Since $\beta_1 > \beta_3$, $\Phi$ is a contraction and thus by the fixed point theorem there indeed exists a unique solution to
\[
    \nu^\rho = \Phi(\nu^\rho).
\]

{\bf Proof of Step 4}. $\Pi^{x, \rho}$ can finally be defined as $\pi^{x, \nu^\rho}$. This family of measures solves the system \eqref{eq:system_for_Pi_rho} by construction. The uniqueness can also be easily checked.

    To finish the proof, we need to show the Lipschitz property of $\Pi^{x, \rho}$ w.r.t. $x$, since the Lipschitz property w.r.t. $\rho$ was already shown in step 1. 
    
    Fix some $x_1, x_2 \in \R^d$, $\rho_1, \rho_2 \in \mathcal{P}_2(\R^d)$. Let $Y_0^1, Y_0^2$ be two $\R^m$--valued random variables with $\Law(Y_0^i) = \Pi^{x_i, \rho_i}$ for $i = 1, 2$ that coupled in an optimal way in the sense of the $\Wass_2$--distance, i.e.:
    \[
        \E|Y_0^1 - Y_0^2|^2 = \Wass_2^2(\Pi^{x_1, \rho_1}, \Pi^{x_2, \rho_2}).
    \]
    Define the following two processes:
    \begin{equation*}
        \dd{Y^i_t} = b(x_i, Y^i_t, \nu^{\rho_i})\dd{t} + \sqrt{2} \dd{W_t},
        \quad i = 1,2,
    \end{equation*}
    where $\nu^{\rho_i} (\dd{y}) = \int \rho_i(\dd{z}) \Pi^{z, \rho_i} (\dd{y})$ and the initial conditions $Y^i_0$ are the random variables fixed before. 
    By what was shown before, $\Pi^{x_i, \rho_i}$ are stationary distributions of $Y^i$, i.e. $\Law(Y_t^i) = \Pi^{x_i, \rho_i}$ for any $t \geq 0$. 
    Note also that the two processes are defined using the same Brownian motion (synchronous coupling). Then, the following inequality holds:
    \begin{equation*}
        \begin{aligned}
            \frac{\dd}{\dd{t}}\E |Y^1_t - Y^2_t|^2 
            &= 2 \E \big\langle Y^1_t - Y^2_t, 
            b(x_1, Y^1_t, \nu^{\rho_1}) - b(x_2, Y^2_t, \nu^{\rho_2}) \big\rangle \\
            & \leq 2 [-\beta_1 \E|Y^1_t - Y^2_t|^2 
            + \beta_2 |x_1 - x_2|^2 + \beta_3\Wass_2^2(\nu^{\rho_1}, \nu^{\rho_2})].
        \end{aligned}
    \end{equation*}
    Applying Grönwall's lemma yields
    \[
    \E|Y_t^1-Y_t^2|^2
    \le
    \e^{-2\beta_1t}\Wass_2^2(\Pi^{x_1,\rho_1},\Pi^{x_2,\rho_2})
    +
    \frac{\beta_2}{\beta_1}|x_1-x_2|^2
    +
    \frac{\beta_3}{\beta_1}
    \Wass_2^2(\nu^{\rho_1},\nu^{\rho_2}).
    \]
    Since the two processes $Y^i$ are stationary, i.e. 
    \(
    \Law(Y_t^i)=\Pi^{x_i,\rho_i}
    \)
    for any $t \geq 0$,  by \eqref{eq:Wass_2_property},
    we obtain
    \[
    \Wass_2^2(\Pi^{x_1,\rho_1},\Pi^{x_2,\rho_2})
    \le
    \e^{-2 \beta_1 t}
    \Wass_2^2(\Pi^{x_1,\rho_1},\Pi^{x_2,\rho_2})
    +
    \frac{\beta_2}{\beta_1}|x_1-x_2|^2
    +
    \frac{\beta_3}{\beta_1}
    \Wass_2^2(\nu^{\rho_1},\nu^{\rho_2}).
    \]
    This inequality holds for any $t \geq 0$; in particular, by letting $t \to \infty$, we get the following:
    \begin{equation}
        \label{eq:aux:Wass_Pi_bound}
    \Wass_2^2(\Pi^{x_1,\rho_1},\Pi^{x_2,\rho_2})
    \le
    \frac{\beta_2}{\beta_1}|x_1-x_2|^2
    +
    \frac{\beta_3}{\beta_1}
    \Wass_2^2(\nu^{\rho_1},\nu^{\rho_2}).
    \end{equation}
    
    It is left to show the control of $\Wass_2(\nu^{\rho_1}, \nu^{\rho_2})$ by $\Wass_2(\rho_1, \rho_2)$, for which we use the definition of the former. To this end, let $\lambda$ be the optimal coupling between $\rho_1$ and $\rho_2$ and $\gamma_{z_1, z_2}$ be the optimal coupling between $\Pi^{z_1, \rho_1}$ and $\Pi^{z_2, \rho_2}$, i.e.:
    \begin{equation*}
        \int |z_1 - z_2|^2 \lambda(\dd{z_1}, \dd{z_2}) = \Wass_2^2 (\rho_1, \rho_2) 
        \quad \text{and} \quad
        \int |y_1 - y_2|^2 \gamma_{z_1, z_2}(\dd{y_1}, \dd{y_2}) = \Wass_2^2 (\Pi^{z_1, \rho_1}, \Pi^{z_2, \rho_2}).
    \end{equation*}
    Note that the following measure on $\R^m \times\R^m$ is therefore a possible coupling (not necessarily optimal) between $\nu^{\rho_1}$ and $\nu^{\rho_2}$:
    \begin{equation*}
        \Gamma(\dd{y_1}, \dd{y_2}) := \int \gamma_{z_1, z_2} (\dd{y_1}, \dd{y_2}) \lambda(\dd{z_1}, \dd{z_2}).
    \end{equation*}
    It is easy to check that the marginals of $\Gamma$ indeed correspond to $\nu^{\rho_1}$ and $\nu^{\rho_2}$. Therefore, using Equation~\eqref{eq:aux:Wass_Pi_bound}, we get:
    \begin{equation*}
    \begin{aligned}
        \Wass_2^2 (\nu^{\rho_1}, \nu^{\rho_2}) &\leq \int |y_1 - y_2|^2 \Gamma(\dd{y_1}, \dd{y_2}) = \iint |y_1 - y_2|^2\gamma_{z_1, z_2} (\dd{y_1}, \dd{y_2}) \lambda(\dd{z_1}, \dd{z_2}) \\
        & \leq \int \Wass_2^2 (\Pi^{z_1, \rho_1}, \Pi^{z_2, \rho_2}) \lambda(\dd{z_1}, \dd{z_2})  \leq \int \bigg[\frac{\beta_2}{\beta_1} |z_1 - z_2|^2 + \frac{\beta_3}{\beta_1} \Wass_2^2(\nu^{\rho_1}, \nu^{\rho_2}) \bigg] \lambda(\dd{z_1}, \dd{z_2}) \\
        &= \frac{\beta_2}{\beta_1}\Wass_2^2(\rho_1, \rho_2) + \frac{\beta_3}{\beta_1} \Wass_2^2(\nu^{\rho_1}, \nu^{\rho_2}).
    \end{aligned}
    \end{equation*}
    Since $\beta_1 > \beta_3$ by Assumption~\ref{assu:gen}, this finally gives us the following upper bound:
    \begin{equation*}
        \Wass_2^2 (\nu^{\rho_1}, \nu^{\rho_2}) \leq \frac{\beta_2}{\beta_1 - \beta_3} \Wass_2^2(\rho_1, \rho_2).
    \end{equation*}
    Plugging it back into Equation~\eqref{eq:aux:Wass_Pi_bound} gives us
    \[
        \Wass_2^2(\Pi^{x_1, \rho_1}, \Pi^{x_2, \rho_2}) \leq \frac{\beta_2}{\beta_1} |x_1 - x_2|^2 + \frac{\beta_2\beta_3}{\beta_1 (\beta_1 - \beta_3)} \Wass_2^2(\rho_1, \rho_2),
    \]
    which proves the Lipschitz continuity by choosing $\widetilde{K}_1 := \frac{\beta_2}{\beta_1}$ and $\widetilde{K}_2 = \frac{\beta_2\beta_3}{\beta_1 (\beta_1 - \beta_3)}$.

\section{Proof of the Main Result}
\label{s:proof_main_result}

We divide the interval $[0, T]$ on $\lceil T/\delta \rceil$ intervals of size at most $\delta$. On each such interval the slow process is expected to not change much and we thus approximate the fast process $Y^{\eps}$ by the process  $\hat{Y}^{\delta, \eps}$ which is successively defined on each time interval $[k \delta , (k + 1) \delta \wedge T)$, where $k = 0, \dots, \lfloor T/\delta \rfloor$, as
\begin{equation*}
\begin{cases}
    \dd{\hat{Y}^{\delta, \eps}_t} & = \frac{1}{\eps} b(X^{\eps}_{k \delta}, \hat{Y}^{\delta, \eps}_t, \hat{\nu}^{\delta, \eps}_t) \dd{t} + \frac{\sqrt{2}}{\sqrt{\eps}} \dd{W_t}, 
    \quad t \in [k \delta, (k + 1)\delta),\\
    \hat{\nu}^{\delta, \eps}_t &= \mathcal{L}(\hat{Y}^{\delta, \eps}_t), \\
    \hat{Y}^{\delta, \eps}_{k \delta} &= Y^\eps_{k \delta} \text{ a.s.}
\end{cases}
\end{equation*}

and let

\begin{equation*}
    \dd{\hat{X}_t^{\delta, \varepsilon}} 
    = v(\hat{X}_t^{\delta, \varepsilon}, \hat{Y}_t^{\delta, \varepsilon})\dd{t}  + \sqrt{2\kappa} \dd{B_t},
    \quad \hat{X}_0^{\delta, \varepsilon} = x_0.
\end{equation*}
We choose $\hat{Y}^{\delta, \eps}$ to be driven by the same Brownian motion as $Y^\eps$. 
Similarly, the process $\hat{X}^{\delta, \varepsilon}$ is coupled to $X^\eps$ using the same Brownian motion $B$. 
Note also that $\hat{Y}^{\delta, \eps}$ coincides with $Y^{\eps}$ at the points of the mesh, while $\hat{X}^{\delta, \varepsilon}$ does not coincide with $X^{\eps}$ at the mesh. 

The proof of the main theorem trivially follows from the following two estimates.
\begin{lemma}
\label{lm:X-X_hat_gen}
    For any $T > 0$, there exists a constant $C_1 > 0$, such that for any $\delta > 0$, we have
    \begin{equation*}
        \sup_{\eps < 1}\sup_{0 \leq t \leq T} \E |X^\eps_t - \hat{X}_t^{\delta, \eps}|^2 \leq 
        \begin{cases}
            C_1 \delta^2, \text{ if } \kappa = 0, \\
            C_1 \delta, \text{ if } \kappa > 0.
        \end{cases}
    \end{equation*}
\end{lemma}

\begin{lemma}
\label{lm:X_hat-X_bar_gen}
    For any $T > 0$, there exists a constant $C_2 > 0$, such that for any $\delta > 0$, we have
    \begin{equation*}
        \lim_{\eps \to 0} \sup_{0 \leq t \leq T} \E |\hat{X}^{\delta, \eps}_t - \bar{X}_t|^2 \leq
         \begin{cases}
            C_2 \delta^2, \text{ if } \kappa = 0, \\
            C_2 \delta, \text{ if } \kappa > 0.
        \end{cases}
    \end{equation*}
\end{lemma}
\begin{rem} 
    Note that the conclusion does not depend on the value of $\kappa$; only the convergence rate in $\delta$ differs: it is of order $\delta^2$ in the absence of Brownian noise and of order $\delta$ otherwise. The constants $C_1$ and $C_2$ may depend on $\kappa$; the case distinction in the estimates is intended only to highlight the different convergence rates in the presence and absence of Brownian noise. In the proofs below, we write estimates of this type in the form 
    \[ 
        C(\delta^2 + \kappa\delta), 
    \] 
    where $C > 0$ may depend on $\kappa$. 
\end{rem}

\begin{proof}[Proof of main Theorem~\ref{thm:main_thm}]
    Fix $T > 0$. By Lemmas~\ref{lm:X-X_hat_gen} and \ref{lm:X_hat-X_bar_gen}, there exist $C_1, C_2 > 0$ such that for any $\delta > 0$
    \begin{equation*}
    \begin{aligned}
        \lim_{\eps \to 0} \sup_{0 \leq t \leq T} &\E |X_t^\eps - \bar{X}_t|^2 
        \leq 
        \lim_{\eps \to 0}\Big[ 
            2\sup_{0 \leq t \leq T} \E |X^\eps_t - \hat{X}_t^{\delta, \eps}|^2 
            +
            2\sup_{0 \leq t \leq T} \E |\hat{X}^{\delta, \eps}_t - \bar{X}_t|^2
            \Big] \\
        & \leq 
        2\sup_{\eps < 1} \sup_{0 \leq t \leq T} \E |X^\eps_t - \hat{X}_t^{\delta, \eps}|^2 
        +
        2\lim_{\eps \to 0} \sup_{0 \leq t \leq T} \E |\hat{X}^{\delta, \eps}_t - \bar{X}_t|^2 
        \leq 2(C_1 + C_2) (\delta^2 + \kappa \delta).
    \end{aligned}
    \end{equation*}
    Since this upper bound holds for any $\delta > 0$, this proves the theorem.
\end{proof}

In the following, we will prove Lemmas~\ref{lm:X-X_hat_gen} and \ref{lm:X_hat-X_bar_gen}. We will start by proving some a priori estimates for the processes involved in the proof.

\subsection{A priori estimates}

\begin{lemma}
\label{lem:X_regulatity_of_path_gen}
    There exists $C^\prime > 0$ such that for any $0 < h < 1$,
    \begin{equation*}
        \sup_{0 \leq t \leq T} \E |X^\eps_{t + h} - X^\eps_t|^2 \leq C^\prime (h^2 + \kappa h).
    \end{equation*}
\end{lemma}

\begin{proof}
\begin{equation*}
\begin{aligned}
        \E |X^\eps_{t + h} - X^\eps_t|^2  
        &= \E \Big| \int_t^{t + h} v(X_s^\eps, Y_s^\eps) \dd{s}  + \sqrt{2\kappa} (B_{t + h} - B_t) \Big|^2 \\
        & \leq 2\E \Big| \int_t^{t + h} v(X_s^\eps, Y_s^\eps) \dd{s}\Big|^2 + 4 \kappa \E |B_{t + h} - B_t|^2\\
        &\leq h \int_t^{t + h} 2K_4(1 + \E|X_s^\eps|^2 + \E|Y_s^\eps|^2) \dd{s} + 4\kappa dh,
\end{aligned}
\end{equation*}
where we used Assumption~\ref{assu:gen}, H\"older's inequality and a simple property of the Brownian motion (recall that $d$ is the dimensionality of the Brownian motions $B$). Using Proposition~\ref{prop:existence_uniqueness}, we obtain
\begin{equation*}
    \E |X^\eps_{t + h} - X^\eps_t|^2 \leq C^\prime (h^2 + \kappa h).
\end{equation*}
\end{proof}

\begin{lemma}
    \label{lem:b_bar_Lipschitz_gen}
    There exist $\widetilde{K}_3, \widetilde{K}_4, \widetilde{K}_5 > 0$ such that for any $x_1, x_2 \in \R^d$, $\rho_1, \rho_2 \in \mathcal{P}_2(\R^d)$,
    \begin{equation*}
        |\bar{v}(x_1, \rho_1) - \bar{v}(x_2, \rho_2) |^2 \leq \widetilde{K}_3 |x_1 - x_2|^2 + \widetilde{K}_4\Wass_2^2(\rho_1, \rho_2)
    \end{equation*}
    and
    \begin{equation*}
        |\bar{v}(x, \rho)|^2 \leq \widetilde{K}_5 (1 + |x|^2 + \Wass_2^2(\rho, \delta_0)).
    \end{equation*}
\end{lemma}
\begin{proof}
    It follows immediately from the definition of $\bar{v}$, Assumption~\ref{assu:gen} and Proposition~\ref{prop:Pi_x_rho_existence}. Indeed, let $\gamma$ be an optimal coupling between $\Pi^{x_1, \rho_1}$ and $\Pi^{x_2, \rho_2}$ and consider
    \begin{equation*}
    \begin{aligned}
        &|\bar{v}(x_1, \rho_1) - \bar{v}(x_2, \rho_2) |^2 = \bigg|\int \big( v(x_1, y_1) - v(x_2, y_2) \big) \gamma(\dd{y_1}, \dd{y_2}) \bigg|^2 \\
        &\leq K_1 |x_1 - x_2|^2 + K_2 \int |y_1 - y_2|^2 \gamma(\dd{y_1}, \dd{y_2}) = K_1 |x_1 - x_2|^2 + K_2 \Wass_2^2(\Pi^{x_1, \rho_1}, \Pi^{x_2, \rho_2})\\
        & \leq (K_1 + K_2\widetilde{K}_1) |x_1 - x_2|^2 + K_2 \widetilde{K}_2 \Wass_2^2(\rho_1, \rho_2),
    \end{aligned}
    \end{equation*}
    which proves the lemma if we define $\widetilde{K}_3 = K_1 + K_2\widetilde{K}_1$ and $\widetilde{K}_4 = K_2 \widetilde{K}_2$.

    The last result follows simply by considering 
    \begin{equation*}
    \begin{aligned}
        |\bar{v}(x, \rho)|^2 &\leq 2|\bar{v}(x, \rho) - \bar{v}(0, \delta_0)|^2 + |\bar{v}(0, \delta_0)|^2 \\
        & \leq 2 \max\{\widetilde{K}_3, \widetilde{K}_4, |\bar{v}(0, \delta_0)|^2\} \big(1 + |x|^2 + \Wass_2^2(\rho, \delta) \big),
    \end{aligned}
    \end{equation*}
    which proves the lemma by taking $\widetilde{K}_5 = 2 \max\{\widetilde{K}_3, \widetilde{K}_4, |\bar{v}(0, \delta_0)|^2\}$.
\end{proof}

\begin{lemma}
    \label{lm:X_bar_bounded_second_moment}
    There exists a constant $\widetilde{C} > 0$ such that
     \begin{equation*}
         \sup_{0 \leq t \leq T } \E |\bar{X}_t|^2 < \widetilde{C},
     \end{equation*}
\end{lemma}
\begin{proof}

    By Lemma~\ref{lem:b_bar_Lipschitz_gen}, and Young's inequality,
    \begin{equation*}
        \begin{aligned}
            \frac{\dd}{\dd{t}} \E|\bar{X}_t|^2 
            &= 2 \E \big\langle \bar{X}_t, \bar{v}(\bar{X}_t, \bar{\rho}_t) \big\rangle 
            + 2\kappa d  
            \leq \E |\bar{X}_t|^2 
            + \E | \bar{v}(\bar{X}_t, \bar{\rho}_t)|^2 
            + 2 \kappa d \\
            & \leq \E |\bar{X}_t|^2 
            + 2\E |\bar{v}(0, \delta_0)|^2 + 4\E |\bar{X}_t|^2 
            + 2\kappa d  \leq C (C + \E |\bar{X}_t|^2),
        \end{aligned}
    \end{equation*}
    where $C$ is a generic constant. The lemma is thus proved by Grönwall's inequality.
\end{proof}

\begin{lemma}
\label{lem:X_bar_regulatity_of_path_gen}
    There exists $K_6 > 0$ such that for any $0 < h < 1$,
    \begin{equation*}
        \sup_{0 \leq t \leq T} \E|\bar{X}_{t + h} - \bar{X}_t|^2 \leq K_6 (h^2 + \kappa h)
    \end{equation*}
    and
    \begin{equation*}
        \sup_{0 \leq t \leq T} \Wass_2^2(\bar{\rho}_{t + h}, \bar{\rho}_t) \leq K_6 (h^2 + \kappa h).\footnote{
            In the above and throughout we should have specified 
            \begin{equation*}
            \sup_{0 \leq t \leq T} \E|\bar{X}_{(t + h) \wedge T} - \bar{X}_t|^2 \leq K_6 (h^2 + \kappa h)
            \quad
            \text{and}
            \quad
            \sup_{0 \leq t \leq T} \Wass_2^2(\bar{\rho}_{(t + h) \wedge T}, \bar{\rho}_t) \leq K_6 (h^2 + \kappa h)
            \end{equation*}
            but we omit such a specification.
        }
    \end{equation*}
\end{lemma}

\begin{proof}
    \begin{equation*}
        \begin{aligned}
            \E|\bar{X}_{t + h} - \bar{X}_t|^2 &= \E \bigg| \int_{t}^{t + h} \bar{v}(\bar{X}_s, \bar{\rho}_s)\dd{s} + \sqrt{2}\kappa (B_{t + h} - B_t) \bigg|^2 \\ 
            &\leq 2\E \bigg| \int_{t}^{t + h} \bar{v}(\bar{X}_s, \bar{\rho}_s) \dd{s} \bigg|^2 + 4 \kappa \E |B_{t + h} - B_t|^2 \\
            & \leq h \int_t^{t + h} 2 \widetilde{K}_5 (1 + \E|\bar{X}_s|^2 + \Wass_2^2(\bar{\rho}_s, \delta_0)) + 4\kappa d h,
        \end{aligned}
    \end{equation*}
    where, in the last inequality, we used Hölder's inequality and Lemma~\ref{lem:b_bar_Lipschitz_gen}. Note that $\Wass_2^2(\bar{\rho}_s, \delta_0) = \E|\bar{X}_s|^2$. Using Lemma~\ref{lm:X_bar_bounded_second_moment}, that gives the uniform control of the second moments, we finally get:
    \[
        \E|\bar{X}_{t + h} - \bar{X}_t|^2 
        \leq 
        K_6 (h^2 + \kappa h)
    \]
    and this upper bound is independent of $t \leq T$. The second result holds trivially, since
    \[
        \Wass_2^2(\bar{\rho}_{t + h}, \bar{\rho}_t) \leq \E|\bar{X}_{t + h} - \bar{X}_t|^2
    \]
    by Equation~\eqref{eq:Wass_2_property}.
\end{proof}

\subsection{Proof of Lemma~\ref{lm:X-X_hat_gen}}

First, we need to show that $Y^\eps$ and $\hat Y^{\delta, \eps}$ are close.
\begin{lemma}
\label{lem:Y-Y_hat_gen}
    There exists $C^{\prime\prime} > 0$ such that
    \begin{equation*}
        \sup_{0 \leq t \leq T} \E |Y^\eps_t - \hat{Y}^{\delta, \eps}_t|^2 \leq C^{\prime\prime}(\delta^{2} + \kappa \delta).
    \end{equation*}
\end{lemma}
\begin{proof}
    For $t \in [k \delta, (k + 1)\delta \wedge T]$, since the two processes are defined with the same Brownian motion, we have:
    \begin{equation*}
        \E |Y^\eps_t - \hat{Y}^{\delta, \eps}_t|^2 
        = 2 \frac{1}{\eps} 
        \int_{k \delta}^t 
        \E \big\langle Y^\eps_s - \hat{Y}^{\delta, \eps}_s, 
         b(X^\eps_s, Y^\eps_s, \nu^\eps_s)
         -  b(X_{k \delta}^\eps, \hat{Y}^{\delta, \eps}_s, \hat{\nu}_s^{\delta, \eps}) \big\rangle 
         \dd{s}. 
    \end{equation*}
    Calling $\varphi^\eps_t := \E |Y^\eps_t - \hat{Y}^{\delta, \eps}_t|^2$ and using Assumption~\ref{assu:gen}, 
    for any $t \in [k \delta, (k + 1)\delta \wedge T]$, we get:
        \begin{equation*}
        \frac{\dd}{\dd{t}}\varphi^\eps_t 
        \leq 
        \frac{1}{\eps} \E 
        \Big[ 
        - \beta_1  |Y^\eps_t - \hat{Y}^{\delta, \eps}_t|^2 
        + \beta_2 |X_t^\eps - X_{k \delta}^\eps|^2 
        + \beta_3 \Wass_2^2(\nu^\eps_s, \hat{\nu}_s^{\delta, \eps}) 
        \Big].
    \end{equation*}
    By Lemma~\ref{lem:X_regulatity_of_path_gen} and Equation~\eqref{eq:Wass_2_property}, we get:
    \begin{equation}
    \label{eq:aux:Y_eps-Y_delta_eps}
        \frac{\dd}{\dd{t}}\varphi^\eps_t 
        \leq 
        -\frac{\beta_1 - \beta_3}{\eps} \varphi^\eps_t  
        + \frac{C^\prime \beta_2}{\eps} (\delta^{2} + \kappa \delta). 
    \end{equation}
    Since $\varphi_{k \delta}^\eps = 0$ for any $k$, Grönwall's inequality gives:
    \begin{equation*}
        \varphi_t^\eps \leq \frac{C^\prime \beta_2}{\beta_1 - \beta_3} (\delta^{2} + \kappa \delta) 
        \big(1 - \e^{-\frac{\beta_1 - \beta_3}{\eps} (t - k \delta)} \big) 
        \leq \frac{C^\prime \beta_2}{\beta_1 - \beta_3} (\delta^{2} + \kappa \delta)
        =: C^{\prime\prime} (\delta^{2} + \kappa \delta),
    \end{equation*}
    uniformly for $t \in [0, T]$, which proves the lemma.
\end{proof}

We can finally proceed to the proof of Lemma~\ref{lm:X-X_hat_gen}.

\begin{proof}[Proof of Lemma~\ref{lm:X-X_hat_gen}]
    Since the two processes $X^\eps$ and $\hat{X}^{\delta, \eps}$ are defined using the same Brownian motion, we have
    \[
        \E |X^\eps_t - \hat{X}_t^{\delta, \eps}|^2 
        = 2 \int_0^t \E \big\langle X^\eps_s - \hat{X}_s^{\delta, \eps}, 
        v(X_s^\eps, Y_s^\eps) -   v(\hat{X}_s^{\delta, \eps}, \hat{Y}_{s}^{\delta, \eps}) \big\rangle \dd{s},
    \]
   By Young's inequality and using Assumption~\ref{assu:gen}, the expression above is upper  bounded by: 
    \[
        (1 + K_1) \int_0^t \E |X^\eps_s - \hat{X}_s^{\delta, \eps}|^2 \dd{s} + K_2 \int_0^t \E|Y_s^\eps - \hat{Y}_{s}^{\delta, \eps}|^2 \dd{s}.
    \]
    By Lemma \ref{lem:Y-Y_hat_gen}, this gives us the following upper bound:
    \begin{equation*}
        \sup_{0 \leq t \leq T} \E |X^\eps_t - \hat{X}_t^{\delta, \eps}|^2
        \leq \int_0^T \sup_{s \leq t} \E |X^\eps_s - \hat{X}_s^{\delta, \eps}|^2 \dd{t} + K_2 C^{\prime \prime} T (\delta^{2} + \kappa \delta).
    \end{equation*}
    By Grönwall's inequality, that finally gives us:
    \begin{equation*}
        \sup_{0 \leq t \leq T} \E |X^\eps_t - \hat{X}_t^{\delta, \eps}|^2 
        \leq
        K_2 C^{\prime \prime} T\e^{T} (\delta^{2} + \kappa \delta),
    \end{equation*}
    which proves the lemma.
\end{proof}

\subsection{Proof of Lemma~\ref{lm:X_hat-X_bar_gen}}

We approximate the process $\hat{Y}^{\delta, \eps}$ by the process $\widetilde{Y}^{\delta, \eps}$ defined below in \eqref{def:Y_tilde}. The difference between \(\widetilde{Y}^{\delta, \eps}\) and \(\hat{Y}^{\delta, \eps}\) is that, in the former, the fixed position of the slow variable takes the value of $\bar{X}_{k \delta}$ instead of $X^\eps_{k \delta}$ and therefore is independent of $\eps$. Consider:
\begin{equation}
\label{def:Y_tilde}
\begin{cases}
    \dd{\widetilde{Y}^{\delta, \eps}_t} & = \frac{1}{\eps} b(\bar{X}_{k \delta}, \widetilde{Y}^{\delta, \eps}_t, \widetilde{\nu}^{\delta, \eps}_t) \dd{t} + \frac{\sqrt{2}}{\sqrt{\eps}} \dd{W_t}, \\
    \widetilde{\nu}^{\delta, \eps}_t &= \mathcal{L}(\widetilde{Y}^{\delta, \eps}_t), \\
    \widetilde{Y}^{\delta, \eps}_{k \delta} &= Y^\eps_{k \delta} \text{ a.s.}
\end{cases}
\end{equation}
Note that, on the time interval $[k \delta, (k + 1)\delta)$, $\widetilde{Y}^{\delta, \eps}$ is substantially a time change of $\mathcal{Y}^{\xi, \zeta}$ defined in \eqref{eq:frozen_process_1_gen}, as we can write:
\begin{equation}
\label{eq:Y_tilde_rescaled_to_Y}
    \widetilde{Y}^{\delta, \eps}_t 
    = 
    \mathcal{Y}^{\bar{X}_{k \delta}, Y^\eps_{k \delta}}_{(t - k\delta)/\eps} 
    \quad 
    \text{a.s. for } t \in [k \delta, (k + 1)\delta \wedge T) ,
\end{equation}
provided we set $\widetilde{W}_{s} = W_{k \delta + s} - W_{k \delta}$ (where $\widetilde{W}$ is the Brownian motion used in the definition of $\mathcal{Y}^{\xi, \zeta}$). The equality is still true if we do not associate the Brownian motions but only in distribution. The condition of independence of $(\xi, \zeta)$ from the Brownian motion (which was part of the definition of the process \eqref{eq:frozen_process_1_gen}) also holds since both $(\bar{X}_{k \delta}, Y^\eps_{k \delta})$, while being dependent on the Brownian motion $W$ up to time $k \delta$, are, in fact, jointly independent of the $\sigma$-algebra generated by $(W_{k \delta + s} - W_{k \delta})_{s \geq 0}$. 

Consider the following lemma that estimates how close $\hat{Y}^{\delta, \eps}$ is to $\widetilde{Y}^{\delta, \eps}$.

\begin{lemma}
\label{lem:Y_hat-Y_tilde_gen}
    There exist constants $C_1, C_2 > 0$ such that
    \begin{equation*}
        \sup_{s \leq t} \E |\hat{Y}^{\delta, \eps}_s - \widetilde{Y}_s^{\delta, \eps}|^2 \leq C_1 (C_2 (\delta^{2} + \kappa \delta) + \sup_{s \leq t} \E |\hat{X}_s^{\delta, \eps} - \bar{X}_s|^2)
    \end{equation*}
\end{lemma}
\begin{proof}
    Since $\hat{Y}^{\delta, \eps}$ and $\widetilde{Y}^{\delta, \eps}$ are driven by the same Brownian motion, the following holds:
    \begin{equation*}
        \E |\hat{Y}^{\delta, \eps}_t - \widetilde{Y}_t^{\delta, \eps}|^2 
        = 
        \frac{2}{\eps} 
        \E 
        \int_0^t \big\langle\hat{Y}^{\delta, \eps}_s - \widetilde{Y}_s^{\delta, \eps}, b(X_{\delta(s)}^\varepsilon, \hat{Y}_s^{\delta, \varepsilon}, \hat \nu_s^{\delta, \varepsilon}) -  b(\bar{X}_{\delta(s)}, \widetilde{Y}_s^{\delta, \eps}, \widetilde{\nu}^{\delta, \eps}_s)\big\rangle \dd{s},
    \end{equation*}
    where $\delta(s) := \lfloor \frac{s}{\delta} \rfloor \delta$ and this notation will be used throughout the section.
    
    Let $\varphi_t^\eps := \E |\hat{Y}^{\delta, \eps}_t - \widetilde{Y}_t^{\delta, \eps}|^2$. Using Assumption~\ref{assu:gen}, we get
    \[
        \frac{\dd}{\dd{t}}\varphi_t^\eps 
        \leq 
        \frac{1}{\eps} 
        \bigg[
        -\beta_1 \varphi_t^\eps 
        + \beta_2 \sup_{s \leq t} \E|X^\eps_s - \bar{X}_s|^2 
        + \beta_3 \Wass_2^2(\Law(\hat{Y}_t^{\delta, \eps}), \Law(\widetilde{Y}_t^{\delta, \eps}))
        \bigg],      
    \]
     Using Equation~\ref{eq:Wass_2_property} and Lemma~\ref{lm:X-X_hat_gen}, we get the following estimate:
    \begin{equation*}
    \begin{aligned}
        \frac{\dd}{\dd{t}}\varphi_t^\eps 
        &\leq 
        \frac{1}{\eps} \bigg[ 
        - (\beta_1 - \beta_3) \varphi_t^\eps 
        + \beta_2 \sup_{s \leq t} \E |X^\eps_s - \bar{X}_s|^2 \bigg]  \\
        &\leq \frac{1}{\eps} \bigg[ 
        - (\beta_1 - \beta_3) \varphi_t^\eps 
        + 2\beta_2 (C_1 (\delta^{2} + \kappa \delta) 
        + \sup_{s \leq t} \E |\hat{X}_s^{\delta, \eps} - \bar{X}_s|^2) \bigg].
    \end{aligned}
    \end{equation*}
    
    By Grönwall's inequality, we finally get
    \[
        \sup_{s \leq t} \E |\hat{Y}^{\delta, \eps}_s - \widetilde{Y}_s^{\delta, \eps}|^2 
        \leq 
        \frac{2\beta_2}{\beta_1 - \beta_3} 
        (C_1 (\delta^{2} + \kappa \delta)
        + \sup_{s \leq t} \E |\hat{X}_s^{\delta, \eps} - \bar{X}_s|^2),
    \]
    which proves the lemma.
\end{proof}

We can now prove Lemma~\ref{lm:X_hat-X_bar_gen}.

\begin{proof}[Proof of Lemma~\ref{lm:X_hat-X_bar_gen}]
    In this proof, $C$ will denote a generic constant, the value of which may change from line to line, but is independent of $\delta$. Recall the notation  $\delta(s) = \lfloor \frac{s}{\delta} \rfloor \delta$. Consider the following inequalities :
    \begin{equation*}
    \begin{aligned}
        \E |\hat{X}^{\delta, \eps}_t - \bar{X}_t|^2 &= \E \Big|\int_0^t v(\hat{X}^{\delta, \eps}_s, \hat{Y}^{\delta, \eps}_s)\dd{s} - \int_0^t \bar{v}(\bar{X}_s, \bar{\rho}_s) \dd{s} \Big|^2 \\
        & \leq 4 \E \Big| \int_0^t v(\hat{X}^{\delta, \eps}_s, \hat{Y}^{\delta, \eps}_s)\dd{s} - \int_0^t v(\bar{X}_s, \widetilde{Y}^{\delta, \eps}_s)\dd{s} \Big|^2 \\
        & \quad + 4 \E \Big|\int_0^t v(\bar{X}_s, \widetilde{Y}^{\delta, \eps}_s)\dd{s} - \int_0^t v(\bar{X}_{\delta(s)}, \widetilde{Y}^{\delta, \eps}_s)\dd{s} \Big|^2 \\
        & \quad + 4 \E \Big|\int_0^t v(\bar{X}_{\delta(s)}, \widetilde{Y}^{\delta, \eps}_s)\dd{s} -  \int_0^t \bar{v}(\bar{X}_{\delta(s)}, \bar{\rho}_{\delta(s)}) \dd{s}  \Big|^2 \\
        & \quad + 4 \E \Big|\int_0^t \bar{v}(\bar{X}_{\delta(s)}, \bar{\rho}_{\delta(s)}) \dd{s} -  \int_0^t \bar{v}(\bar{X}_s, \bar{\rho}_s) \dd{s} \Big|^2 \\
        & =: I_1(t) + I_2(t) + I_3(t) + I_4(t).
    \end{aligned}
    \end{equation*}
    Let us consider $I_1$, $I_2$, $I_3$, and $I_4$ separately. By Assumption~\ref{assu:gen} and Lemma~\ref{lem:Y_hat-Y_tilde_gen}:
    \begin{equation*}
    \begin{aligned}
        I_1(t) &\leq 4T \Big(K_1 \int_0^t \E|\hat{X}_s^{\delta, \eps} - \bar{X}_s|^2 \dd{s} + K_2 \int_0^t \E |\hat{Y}^{\delta, \eps}_s - \widetilde{Y}^{\delta, \eps}_s|^2\dd{s} \Big) \\
        & \leq C \int_0^t \sup_{u \leq s}  \E|\hat{X}_u^{\delta, \eps} - \bar{X}_u|^2 \dd{s} + C (\delta^{2} + \kappa \delta).
    \end{aligned}
    \end{equation*}
    The bound for $I_2$ is a straightforward consequence of Lemma~\ref{lem:X_bar_regulatity_of_path_gen}:
    \begin{equation*}
        \sup_{0 \leq t \leq T} I_2(t) \leq C (\delta^{2} + \kappa \delta).
    \end{equation*}
    For $I_3$, note that we can represent it as:
    \begin{equation*}
        \frac{1}{4} I_3(t) = \E \Big| \sum_{k} \int_{k\delta}^{(k + 1)\delta} v(\bar{X}_{k \delta}, \widetilde{Y}_s^{\delta, \eps}) \dd{s} - \sum_{k} \bar{v}(\bar{X}_{k \delta}, \bar{\rho}_{k \delta}) \delta \Big|^2.
    \end{equation*}
    
    The equality \eqref{eq:Y_tilde_rescaled_to_Y} allows us to rewrite $I_3$ as:
    \begin{equation*}
    \begin{aligned}
        \frac{1}{4} I_3(t) &= \E \Big| \sum_{k = 0}^{\delta(t)} \int_{k\delta}^{(k + 1)\delta} v(\bar{X}_{k \delta}, \mathcal{Y}^{\bar{X}_{k \delta}, Y^\eps_{k \delta}}_{\frac{s - k\delta}{\eps}}) \dd{s} - \sum_{k = 0}^{\delta(t)} \bar{v}(\bar{X}_{k \delta}, \bar{\rho}_{k \delta}) \delta \Big|^2 \\      
        &=\E \Big|\eps \sum_{k = 0}^{\delta(t)} \int_0^{\delta/\eps} v(\bar{X}_{k \delta}, \mathcal{Y}_s^{\bar{X}_{k \delta}, Y_{k\delta}^{\eps}}) \dd{s} - \sum_{k = 0}^{\delta(t)} \bar{v}(\bar{X}_{k \delta}, \bar{\rho}_{k \delta}) \delta \Big|^2 \\
        &=\E \Big| \delta \frac{\eps}{\delta} \sum_{k = 0}^{\delta(t)} \int_0^{\delta/\eps} v(\bar{X}_{k \delta}, \mathcal{Y}_s^{\bar{X}_{k \delta}, Y_{k\delta}^{\eps}}) \dd{s} - \sum_{k = 0}^{\delta(t)} \bar{v}(\bar{X}_{k \delta}, \bar{\rho}_{k \delta}) \delta \Big|^2.
    \end{aligned}
    \end{equation*}

    The last transformation allows us to represent the expression above as an $L^2(\Omega; \R)$ norm of a sum of $(\delta(t) + 1)$ random variables, as follows: set
    \[
        A_{k, \eps} 
        := 
        \frac{\eps}{\delta}
        \int_0^{\delta /\eps}
        v(\bar{X}_{k \delta}, \mathcal{Y}_s^{\bar{X}_{k \delta}, Y_{k\delta}^{\eps}}) \dd{s}
        - 
        \bar{v}(\bar{X}_{k \delta}, \bar{\rho}_{k \delta}).
    \]
    Then, by triangular inequality, we have
    \begin{equation}
    \label{eq:I_3_sum_expression}
        \frac{1}{4}
        \sup_{0 \leq t \leq T} I_3 (t) 
        =
        \sup_{0 \leq t \leq T}
        \Biggl| 
            \delta 
            \sum_{k = 0}^{\delta(t)}
            A_{k, \eps} 
        \Biggr|^2_{L^2(\Omega; \R)}
        \!\! \leq
        \delta^2 (\delta(T) + 1) 
        \sum_{k = 0}^{\delta(T)}
        \big| 
            A_{k, \eps} 
        \big|^2_{L^2(\Omega; \R)}        . 
    \end{equation}
    
    We can now use Proposition~\ref{prop:ergodicity_of_Y} to show that the right-hand side of the above tends to zero. 
    Indeed, by Proposition~\ref{prop:existence_uniqueness}, the second moments of the measures $(\Law(Y_{k\delta}^{\eps}))_{k \leq \lfloor T/\delta \rfloor} \subset \mathcal{K}_{\widetilde{C}}$ are uniformly bounded by a constant $\widetilde{C}$ (which depends on $T$, but not on $\eps$ or $\delta$). Therefore, we can use Proposition~\ref{prop:ergodicity_of_Y} for each individual random variable in the sum \eqref{eq:I_3_sum_expression}, which provides the following convergence
    \begin{equation*}
        \frac{\eps}{\delta}\int_0^{\delta/\eps} v(\bar{X}_{k \delta}, \mathcal{Y}_s^{\bar{X}_{k \delta}, Y_{k\delta}^{\eps}}) \dd{s} \xrightarrow[\eps \to 0]{L^2(\Omega; \R)} \bar{v}(\bar{X}_{k \delta}, \bar{\rho}_{k \delta}),
    \end{equation*}
    and, therefore,
    \[
        \big| 
            A_{k, \eps} 
        \big|^2_{L^2(\Omega; \R)} \xrightarrow[\eps \to 0]{} 0.
    \]
    for any $k = 0, \dots, \delta(T)$. It follows that
    \begin{equation*}
        \sup_{0 \leq t \leq T} I_3(t) \xrightarrow[\eps \to 0]{} 0.
    \end{equation*}

    Let us finally consider $I_4$. By Lemma~\ref{lem:b_bar_Lipschitz_gen}, $\bar{v}$ is Lipschitz continuous.
    Together with Lemma~\ref{lem:X_bar_regulatity_of_path_gen}, this gives us the following estimate:
    \begin{equation*}
        \frac{1}{4}
        \sup_{0 \leq t \leq T}
        I_4(t) 
        \leq T \int_0^T \big(K_1\E|\bar{X}_{\delta(s)} - \bar{X}_s|^2 + K_2\Wass_2^2(\bar{\rho}_{\delta(s)}, \bar{\rho}_s)\big) \dd{s} \leq 2 (K_1 + K_2) K_6 T^2 (\delta^{2} + \kappa \delta).  
    \end{equation*}
    
    Combining the estimates for $I_1, I_2, I_3$, and $I_4$, we finally have the following:
    \[
         \sup_{0 \leq t \leq T}
         \E |\hat{X}^{\delta, \eps}_t - \bar{X}_t|^2 \leq C \int_0^T \sup_{u \leq s}  \E|\hat{X}_u^{\delta, \eps} - \bar{X}_u|^2 \dd{s} + C(\delta^{2} + \kappa \delta)
         + 
         \sup_{0 \leq t \leq T} I_3(t). 
    \]
    By Grönwall's inequality that gives:
    \[
        \sup_{0 \leq t \leq T}
        \E |\hat{X}^{\delta, \eps}_t - \bar{X}_t|^2 
        \leq \big(C (\delta^{2} + \kappa \delta) + \sup_{0 \leq t \leq T}I_3(t)\big) \e^{CT}
        \xrightarrow[\eps \to 0]{} C (\delta^{2} + \kappa \delta).
    \]
    
    That proves the lemma.

\end{proof}

\begin{rem}[Comparison with \cite{XLLM21}]\label{note:explnation error}
The proof of the averaging principle in \cite{XLLM21} is based on a
time-discretisation argument similar to the one used in this section. Roughly speaking,
\cite[Lemma~5.4]{XLLM21} plays the role of
Lemma~\ref{lm:X-X_hat_gen}, whereas \cite[Lemma~6.4]{XLLM21}
corresponds to Lemma~\ref{lm:X_hat-X_bar_gen}.

The final part of the argument is also similar in both proofs. On each
interval $[k\delta,(k+1)\delta)$, the corresponding frozen fast process
is identified, after rescaling time, with a process of the form
$\mathcal Y^{\xi,\zeta}$. To be more precise,  in \cite{XLLM21}, the frozen process is identified as being a rescaling of 
$\hat Y$, 
\[
    \mathcal Y^{X^\eps_{k\delta},Y^\eps_{k\delta}},
\]
whereas in our proof we take the frozen process to be a rescaling of 
$\widetilde Y^{\delta, \eps}$, leading to
\[
    \mathcal Y^{\bar X_{k\delta},Y^\eps_{k\delta}}.
\]

The introduction of $\widetilde Y^{\delta, \eps}$ is merely a
technical difference and is not the essential point of the comparison.
The important issue lies instead in how the process
$\mathcal Y^{X^\eps_{k\delta},Y^\eps_{k\delta}}$ is treated in
\cite{XLLM21}. In the proof of \cite[Lemma~6.2]{XLLM21}, a long-time convergence estimate is
first established for a family of frozen processes
$\mathcal Y^{x,y}$ indexed by deterministic values
$(x,y)\in\R^d\times\R^m$. The random variables
$(X^\eps_{k\delta},Y^\eps_{k\delta})$ are then substituted for
$(x,y)$ in this estimate; see the way the quantity named there $\mathcal{J}_k^{\epsilon}$ is treated, especially \cite[p.~S860, line~4]{XLLM21}.
In other words, the map
\[
    (x,y)\mapsto\mathcal Y^{x,y}
\]
is evaluated pointwise at the random pair
$(X^\eps_{k\delta},Y^\eps_{k\delta})$.

As explained in Note~\ref{rem:Y_xi_zeta_parameters},
unlike in the linear case,
the process obtained by the pointwise substitution
\[
    \omega\mapsto
    \mathcal Y^{
        X^\eps_{k\delta}(\omega),
        Y^\eps_{k\delta}(\omega)
    }(\omega)
\]
is not, in general, the same as the McKean--Vlasov process
$\mathcal Y^{X^\eps_{k\delta},Y^\eps_{k\delta}}$ defined directly
with the random parameters
$(X^\eps_{k\delta},Y^\eps_{k\delta})$ through
\eqref{eq:frozen_process_1_gen}. The latter is the actual
time-rescaled version of $\hat Y$ on
$[k\delta,(k+1)\delta]$, and its long-time behaviour depend on the law induced by the first random parameter. Put differently, in \cite{XLLM21} the process $\mathcal Y^{\xi, \zeta}$ is not viewed as a McKean-Vlasov evolution with random drift, but the value of the drift is fixed. As we have argued the law of the variable $\xi$ will have an effect on the long time properties of the frozen process. Without taking this into account, the fast process is treated as if it was a Markov process. 

This substitution seems to lead to an identification of the effective drift 
in \cite{XLLM21} that is, in general, incorrect.
\end{rem}

\appendix

\section{Miscellaneous facts}\label{appendix A}

\subsection{Fokker-Planck for conditional law}\label{appendix: conditional PDE}
The fact that the law  $\varsigma_t^{x, \rho}$ of the process $\mathcal Y^{\xi, \zeta} \vert \xi=x $  satisfies equation 
\eqref{eqn: evol conditional law intro}
can be seen in two ways. Firstly, by looking at equation \eqref {eq:def:mu_xi_zeta} for the joint law and substituting, in that equation, the decomposition $\boldsymbol{\mu}^{\xi, \zeta}(x,y)=\rho(x) \varsigma^{x, \rho}(y)$ it is clear that $\varsigma^{x, \rho}$ solves \eqref{eqn: evol conditional law intro} for $\rho$-almost every $x$.
Alternatively, one can write the It\^{o} formula applied to the process $\mathcal{Y}^{\xi}$ in \eqref{eqn: our frozin proces intro} and then take conditional expectation with respect to $\xi=x$. 

Indeed, let us set $d=m=1$ just to shorten notation; then,  for any $\varphi\in C_c^2(\mathbb{R})$, we have 
    \[
       \dd \varphi(\mathcal{Y}_t^{\xi}) 
        = 
         \varphi'(\mathcal{Y}_t^{\xi})  b(\xi, \mathcal Y_t^{\xi}, \bnu_t^{\xi})
        +  \varphi''(\mathcal{Y}_t^{\xi}) 
        + \sqrt{2}   \varphi(\mathcal{Y}_s^{\xi}) \dd{W_t}. 
    \]
    First take expectation on both sides with respect to the joint law  $\boldsymbol{\mu}^{\xi, \zeta}$, i.e. write 
    \begin{align*}
    \dd\mathbb E_{\boldsymbol{\mu}} \left[\varphi(\mathcal{Y}_t^{\xi}) \right] & = \dd \int_{\R} \rho(\dd x) \int_{\R} \dd y \,  \varphi(y) \varsigma_t^{x,\rho}(y)\\ & = \int_{\R} \rho(\dd x) \int_{\R} \dd y \, \varphi'(y) b(x,y, \bnu_t^{\xi}) \varsigma_t^{x,\rho}(y) +\int_{\R} \rho(\dd x) \int_{\R} \dd y \, \varphi''(y) \varsigma_t^{x,\rho}(y) \,.
    \end{align*}
From the above the weak equation for $\varsigma_t^{x,\rho}(y)$ can be read off, and one can see that such an equation can only hold for $\rho$-almost every $x$.
\subsection{Asymptotics of the law}

In Section~\ref{s:heuristics}, we heuristically derived the effective behavior of $X^\eps$ as $\eps \to 0$ by working directly with the system of SDEs \eqref{eq:def:McKean-Vlasov_full} and using the time-scale separation at the level of trajectories. Here, instead, we derive the effective equation by examining the asymptotic behaviour of the law of the pair $(X^\varepsilon,Y^\varepsilon)$.
The argument presented here is only informal and does not correspond to the one used in the proof, yet it provides a different perspective on the main result.

We begin with the case in which the fast process does not depend on its own law; i.e. the classical averaging setting \eqref{eq:multiscale system intro}. 
Recall that $\eta^\eps_t = \mathcal L(X_t^\varepsilon,Y_t^\varepsilon)$ denotes the law of the pair and it solves the following equation
\[
    \partial_t \eta_t^\varepsilon(x, y)
    =
    -\nabla_x \cdot \left[
    v(x,y)\eta_t^\varepsilon(x,y)
    \right] 
    + \kappa \Delta_{x} \eta_t^\varepsilon(x, y)
    -\frac{1}{\varepsilon}
    \nabla_y \cdot
    \!\left[
    b(x,y)\eta_t^\varepsilon(x,y)
    \right]
    + \frac{1}{\varepsilon} \Delta_{y} \eta_t^\varepsilon(x,y),
\]
that is,
\begin{equation}
\label{eq:standard_law_evolution}
    \partial_t \eta_t^\eps
    = \mathcal L^\varepsilon \eta_t^\eps = \mathcal L_{S}\eta_t^\eps
    +
    \frac{1}{\varepsilon}\mathcal L_{F}\eta_t^\eps,
\end{equation}
where $\mathcal L_{S}$ and $\mathcal L_{F}$ are two linear differential operators corresponding to the slow and the fast components of the dynamics:
\begin{equation*}
    (\mathcal L_{S}\eta)(x,y)
    :=
    -\nabla_x \cdot \big[v(x,y)\eta(x,y)\big]
    +
    \kappa \Delta_{x} \eta(x,y), 
    \;\; 
    (\mathcal L_{F}\eta)(x,y)
    :=
    -\nabla_y \cdot \big[ b(x,y)\eta(x,y)\big]
    +
    \Delta_{y} \eta(x,y).
\end{equation*}
Note that $\mathcal{L}_S$ and $\mathcal{L}_F$ can be seen as operators on $C(\R^d)$ only and $C(\R^m)$ only respectively, while the respective other variable can be seen as a parameter.

Recall that
\[
    \rho_t^\varepsilon(\dd x)
    =
    \int_{\R^m}\eta_t^\varepsilon(\dd x,\dd y)
\]
denotes the marginal law of the slow component. We can decompose
the joint law as
\[
    \eta_t^\varepsilon(\dd x,\dd y)
    =
    \rho_t^\varepsilon(\dd x)\,
    \nu_t^\varepsilon(\dd y|x),
\]
where
\[
    \nu_t^\varepsilon(\dd y|x)
    =
    \Law(Y_t^\varepsilon\,|\,X_t^\varepsilon=x)
\]
is the conditional law of the fast component given the value of the slow variable.

Formally multiplying Equation~\eqref{eq:standard_law_evolution} by
\(\varepsilon\), we see that, if the limit
\(\displaystyle\bar\eta_t:=\lim_{\varepsilon\to0}\eta_t^\varepsilon\) exists, it should satisfy
\[
    \mathcal L_F\bar\eta_t =0.
\]

Since the fast component evolves on a much shorter time scale, it is natural to expect that, for each value of the slow variable $x$, the conditional law rapidly approaches an equilibrium (potentially) associated with that value of $x$. In particular, in the limit $\varepsilon \to 0$, the conditional law should not possess its own time dependence; any variation in time should occur only through the evolution of the slow variable. This motivates the following {\it ansatz} in the case of classical averaging:
\begin{equation}
    \label{eq:ansatz_standard}
    \bar\eta_t(\dd{x}, \dd{y})
    =
    \bar\rho_t(\dd{x})\bar\nu(\dd{y}|x),    
\end{equation}
where the conditional law \(\bar\nu(\cdot|x)\) is independent of \(t\).
Substituting this decomposition into the stationarity condition
\(\mathcal L_F \bar\eta_t=0\), we obtain
\[
    \bar\rho_t\,\mathcal L_F\bar\nu =0.
\]
Hence, for \(\bar\rho_t\)-a.e. \(x\),
\[
    \mathcal L_F\bar\nu =0.
\]
Under the standard assumptions of averaging, this equation admits a unique
solution, denoted by \(\pi^x\), which is precisely the invariant measure
of the frozen process corresponding to the parameter \(x\). 

This heuristics can be justified in view of the recent work \cite{xie2026uniform}, which has looked at the limit of the joint law $\eta^{\eps}_t$ and proved that, as $\eps \rightarrow 0$, it converges to the product $\mathcal L_{\bar {X}_t}(\dd{x}) \pi^x(\dd{y})$.

\begin{rem}
The only additional heuristic assumption made above is that, in the limit $\eps \to 0$, the conditional law of the fast variable given the slow variable becomes independent of time. 
Under this assumption, the stationarity condition $\mathcal L_F \bar\eta_t = 0$ immediately implies that this limit conditional law $\bar{\nu}$ must be an invariant measure of the frozen fast dynamics. Under the uniqueness assumption, it is therefore given by $\pi^x$. In particular, the above argument directly explains why the averaging procedure involves integration with respect to $\pi^x$.

This is not, however, the usual way in which averaging is presented. 
As explained in the introduction, the usual argument starts from the frozen fast process $Y^{(x)}$ and its invariant measure $\pi^x$, whereas the argument above starts from the evolution equation for the joint law and identifies the limit conditional law with $\pi^x$ by assuming that the former is independent of time.

The reason why the latter viewpoint is useful to derive the effective equation for the system \eqref{eq:def:McKean-Vlasov_full} is the following. When the fast process is nonlinear, constructing a frozen equation apriori is not straightforward. As we will see, the law-based viewpoint allows us instead to modify the heuristic assumption that the conditional law becomes independent of time in the limit $\eps \to 0$. This modified ansatz will lead directly to the appropriate limit equation in the nonlinear setting.\footnote{
    It is important to stress that the limit considered here is only the limit $\eps\to0$, with $t$ kept fixed. The distinction between the limits $\eps\to0$ and $t\to\infty$ is important in view of the discussion in \cite{duong2025coarse}, where the authors distinguish, in particular, between averaging and projection methods based on conditional expectations. In the latter approach, the effective dynamics is obtained by projecting onto the equilibrium conditional law of the fast variable given the slow one. These two procedures need not lead to the same effective dynamics, even in the presence of time-scale separation.

    Our approach also bears some resemblance to the Gy\"{o}ngy the method \cite{duong2025coarse}, where the dimension reduction is achieved by integrating the drift of the slow process with respect to the conditional law of the fast process given the slow variable. More precisely, in the notation used here, integrating the equation for $\eta^\eps$ with respect to $y$ after writing
    $$
        \eta_t^\varepsilon(\dd x,\dd y)
        =
        \rho_t^\varepsilon(\dd x)\,
        \nu_t^\varepsilon(\dd y|x)
    $$
    leads to an effective drift obtained by averaging $v(x,y)$ with respect to the conditional law $\nu_t^\varepsilon(\dd y|x)$.
    This is precisely the type of conditional expectation appearing in the Gy\"{o}ngy construction. However, the Gy\"{o}ngy method itself does not prescribe a particular behaviour of the conditional law as $\varepsilon\to0$: the conditional law is, in general, time-dependent.
    This is where our heuristic assumption enters: we assume that, for fixed $t$, the conditional law becomes independent of time in the limit $\eps \to 0$, which allows to pass to the correct limit when $\eps \to 0$.
}
\end{rem}

We can finally describe a closed effective equation for the limit $\displaystyle \bar{\rho}_t = \lim_{\eps \to 0} \rho^\eps_t$. Plugging the decomposition 
\(
    \bar\eta_t(\dd{x}, \dd{y})
    =
    \bar\rho_t(\dd{x}) \pi^x(\dd{y}),
\) 
back in \eqref{eq:standard_law_evolution}, integrating w.r.t. $y$ and passing to the limit $\eps \to 0$, we get
\[
    \partial_t \bar{\rho}_t (x) = 
    -\nabla_x\cdot\left[
    \int v(x,y) \pi^x(\dd{y}) \bar{\rho}_t(x) 
    \right] 
    + \kappa\Delta_{x} \bar{\rho}_t(x),
\]
which corresponds to the law of the process
\[
    \dd{\bar{X}_t} = \bar{v}(\bar{X}_t) \dd{t} + \sqrt{2} \kappa \dd{B_t},
\]
where 
\( \bar{v}(x) = \int_{\R^m} v(x,y)\,\pi^x(\dd y), \)
which is precisely the effective equation obtained in the classical averaging theory for diffusion processes.

In case of the process \eqref{eq:def:McKean-Vlasov_full}, the law instead solves the following equation
\begin{equation*}
\begin{cases}
    \partial_t \eta_t^\varepsilon(x, y)
    =
    -\nabla_x\cdot\left[
    v(x,y)\eta_t^\varepsilon(x,y)
    \right] 
    + \kappa \Delta_{x} \eta_t^\varepsilon(x, y)
    -\frac{1}{\varepsilon}\nabla_y\cdot
    \left[
    b(x,y,\nu_t^\varepsilon)\eta_t^\varepsilon(x,y)
    \right]
    + \frac{1}{\varepsilon} \Delta_{y}\eta_t^\varepsilon(x,y), \\    
    \nu_t^\eps = \int \eta_t^\varepsilon(\dd{z}, y),
\end{cases}
\end{equation*}
which can also be written using the operator notation:
\begin{equation}
\label{eq:semigroup_McKean}
    \partial_t \eta_t^\eps
    = \mathcal L^\varepsilon\!\left[ \eta_t^\eps\right],
\end{equation}
where $\mathcal L^\varepsilon$ is, this time, {\it a non-linear operator} that splits into two parts (slow and fast one):
\[
    \mathcal L^\varepsilon[\eta] =
    \mathcal L_{S}\eta
    +
    \frac{1}{\varepsilon}\mathcal L_{F}[\eta].
\]
The slow part is linear and reads
\[
    (\mathcal L_{S}\eta)(x,y)
    =
    -\nabla_x \cdot \big[v(x,y)\eta(x,y)\big]
    +
    \kappa \Delta_{x} \eta(x,y),
\]
while the fast part is non-linear and is defined as:
\[
    (\mathcal L_{F}[\eta])(x,y)
    =
    -\nabla_y \cdot 
    \big[b(x,y, \int \eta(z, \cdot) \dd{z})\eta(x,y)\big]
    +
    \Delta_{y} \eta(x,y).
\]

As in the standard case, if the limit $\ds \bar{\eta}_t := \lim_{\eps \to 0} \eta^\eps_t$ exists, it should necessarily satisfy:
\[
    \mathcal{L}_F[\bar{\eta}_t] = 0.
\]
Nevertheless, because of the lack of linearity of the fast part, the ansatz \eqref{eq:ansatz_standard} is not expected to hold. The alternative ansatz would be that the conditional distribution may depend on time but only {\it as a function of} $\bar{\rho}_t$ (the law of the effective process):
\begin{equation*}
\label{eq:aux:mu_bar_decomposition}
    \bar{\eta}_t(x, y) = \bar{\rho}_t(x) \bar{\nu}^{\bar{\rho}_t}(y|x).
\end{equation*}
If we now plug this decomposition back into the stationary equation, we have:
\[
    \mathcal{L}_F[\bar{\eta}_t] = \bar{\rho}_t \mathcal{L}_F[\bar{\nu}^{\bar{\rho}_t}(y|x)] = 0,
\]
which means that, for $\bar{\rho}_t$--a.e. $x$, $\bar{\nu}^{\bar{\rho}_t}(y|x)$ solves the following equation:
\[
    -\nabla_y \cdot \bigg[b \big(x,y, \int \bar{\nu}^{\bar{\rho}_t}(\cdot|z) \dd{z} \big) \,\bar{\nu}^{\bar{\rho}_t}(y|x)\bigg]
    +
    \Delta_{y} \bar{\nu}^{\bar{\rho}_t}(y|x) = 0.    
\]
Therefore, the derivation of the effective equation boils down to the analysis of system \eqref{eqn: Pi intro}. 
This is exactly the measure that is used in the definition of the effective process $\bar{X}$. 
As per Section~\ref{s:Assumptions_main_result}, under our assumptions, there is a unique solution to this system of equations and therefore the conditional law $\bar{\nu}$ can be identified as:
\begin{equation*}
    \bar{\nu}^{\bar{\rho}_t}(y|x) = \Pi^{x, \bar{\rho}_t}(y) 
    \quad \text{for } \bar{\rho}_t \text{--a.e. } x.
\end{equation*}

Plugging the decomposition \eqref{eq:aux:mu_bar_decomposition} and the expression for the conditional law back into \eqref{eq:semigroup_McKean}, integrating both sides in $y$ and taking the limit $\eps \to 0$, we get the following effective equation for $\bar{\rho}$:
\[
    \partial_t \bar{\rho}_t (x) 
    = 
    -\nabla_x\cdot\left[
    \int v(x,y) \Pi^{x, \bar{\rho}_t}(\dd{y}) \bar{\rho}_t(x) 
    \right] 
    + \kappa \Delta_{x} \bar{\rho}_t(x),
\]
which is exactly the Fokker-Planck equation corresponding to the law of the effective process \eqref{eq:X_bar_rho_bar}. Note that, in the case of $\kappa = 0$ (no proper Brownian motion for the slow component), $\bar{\rho}_t = \delta_{\bar{X}_t}$, where $\delta$ stands for the Dirac measure. In this case, just as in the standard averaging, the dependence boils down to only dependence on the position.

\begin{small}
\bibliographystyle{plain}
\bibliography{biblioKramer.bib}

@article{xiang2026averaging,
  title={Averaging principles for nonautonomous multiscale McKean-Vlasov stochastic systems},
  author={Xiang, Jie and Qiao, Huijie},
  journal={arXiv preprint arXiv:2606.12820},
  year={2026}
}

@article{xie2026uniform,
  title={Uniform-in-time diffusion approximations for multiscale stochastic systems},
  author={Xie, Longjie and Zhang, Xicheng},
  journal={arXiv preprint arXiv:2604.00692},
  year={2026}
}

@article{pardoux2001poisson,
  title={On the Poisson equation and diffusion approximation. I},
  author={Pardoux, {\'E}tienne and Veretennikov, A Yu},
  journal={Annals of probability},
  pages={1061--1085},
  year={2001},
  publisher={JSTOR}
}

@inproceedings{crisan2026poisson,
  title={Poisson equations with locally-Lipschitz coefficients and uniform in time averaging for stochastic differential equations via strong exponential stability},
  author={Crisan, Dan and Dobson, Paul and Goddard, Benjamin D and Ottobre, Michela and Souttar, Iain},
  booktitle={Annales de l'Institut Henri Poincare (B) Probabilites et statistiques},
  volume={62},
  number={1},
  pages={667--721},
  year={2026},
  organization={Institut Henri Poincar{\'e}}
}

@article{hong2025averaging,
  title={Averaging principle for multiscale stochastic PDEs with fully local monotone coefficients},
  author={Hong, Wei and Liu, Wei and Yang, Luhan and Yang, Shiyuan},
  journal={Discrete and Continuous Dynamical Systems-S},
  volume={18},
  number={10},
  pages={2791--2809},
  year={2025},
  publisher={Discrete and Continuous Dynamical Systems-S}
}

@article{cerrai2017averaging,
  title={Averaging principle for nonautonomous slow-fast systems of stochastic reaction-diffusion equations: the almost periodic case},
  author={Cerrai, Sandra and Lunardi, Alessandra},
  journal={SIAM Journal on Mathematical Analysis},
  volume={49},
  number={4},
  pages={2843--2884},
  year={2017},
  publisher={SIAM}
}

@article{pardoux2003poisson,
  title={On Poisson equation and diffusion approximation 2},
  author={Pardoux, E and Veretennikov, A Yu},
  journal={The Annals of Probability},
  volume={31},
  number={3},
  pages={1166--1192},
  year={2003},
  publisher={Institute of Mathematical Statistics}
}

@article{pardoux2005poisson,
  title={On the Poisson equation and diffusion approximation 3},
  author={Pardoux, Etienne and Veretennikov, A Yu},
  year={2005}
}

@book{ethier2009markov,
  title={Markov processes: characterization and convergence},
  author={Ethier, Stewart N and Kurtz, Thomas G},
  year={2009},
  publisher={John Wiley \& Sons}
}

@preamble{
   "\def\cprime{$'$} "
}

@book {FW,
    AUTHOR = {Freidlin, M. I. and Wentzell, A. D.},
     TITLE = {Random perturbations of dynamical systems},
    SERIES = {Grundlehren der Mathematischen Wissenschaften [Fundamental
              Principles of Mathematical Sciences]},
    VOLUME = {260},
      NOTE = {Translated from the Russian by Joseph Sz{\"u}cs},
 PUBLISHER = {Springer-Verlag},
   ADDRESS = {New York},
      YEAR = {1984},
     PAGES = {viii+326},
      ISBN = {0-387-90858-7},
   MRCLASS = {60H10 (60J60 70L05)},
  MRNUMBER = {MR722136 (85a:60064)},
}

@book{PS08,
  title={Multiscale methods: averaging and homogenization},
  author={Pavliotis, Grigoris and Stuart, Andrew},
  year={2008},
  publisher={Springer Science \& Business Media}
}

@article{XLLM21,
  author  = {Xu, Jie and Liu, Juanfang and Liu, Jicheng and Miao, Yu},
  title   = {{Strong Averaging Principle for Two-Time-Scale Stochastic
             {McKean--Vlasov} Equations}},
  journal = {Applied Mathematics \& Optimization},
  year    = {2021},
  volume  = {84},
  number  = {Suppl 1},
  pages   = {837--867},
  doi     = {10.1007/s00245-021-09787-3}
}

@book{weinan2011principles,
  title={Principles of multiscale modeling},
  author={Weinan, E},
  volume={1534},
  year={2011},
  publisher={Cambridge University Press, Cambridge}
}

@article{khas1968averaging,
  title={The averaging principle for stochastic differential equations},
  author={Khas'minskii, Rafail Z.},
  journal={Problemy Peredachi Informatsii},
  volume={4},
  number={2},
  pages={86--87},
  year={1968},
  publisher={Russian Academy of Sciences, Branch of Informatics, Computer Equipment and~…}
}

@article{CZ16,
author = {Ren{\'e} Carmona and Xiuneng Zhu},
title = {{A probabilistic approach to mean field games with major and minor players}},
volume = {26},
journal = {The Annals of Applied Probability},
number = {3},
publisher = {Institute of Mathematical Statistics},
pages = {1535 -- 1580},
year = {2016},
doi = {10.1214/15-AAP1125},
URL = {https://doi.org/10.1214/15-AAP1125}
}

@article{legoll2010effective,
  title={Effective dynamics using conditional expectations},
  author={Legoll, Fr{\'e}d{\'e}ric and Lelievre, Tony},
  journal={Nonlinearity},
  volume={23},
  number={9},
  pages={2131--2163},
  year={2010}
}

@book{villani2009optimal,
  title={Optimal transport: old and new},
  author={Villani, C{\'e}dric and others},
  volume={338},
  year={2009},
  publisher={Springer}
}

@article{meyn1993stability,
  title={{Stability of Markovian processes III: Foster--Lyapunov criteria for continuous-time processes}},
  author={Meyn, Sean P and Tweedie, Richard L},
  journal={Advances in Applied Probability},
  volume={25},
  number={3},
  pages={518--548},
  year={1993},
  publisher={Cambridge University Press}
}

@article{chaintron2022propagation,
  title={Propagation of chaos: a review of models, methods and applications. I. Models and methods},
  author={Chaintron, Louis-Pierre and Diez, Antoine},
  journal={arXiv preprint arXiv:2203.00446},
  year={2022}
}

@article{LietalPoisson,
author = {Li, Yun and Wu, Fuke and Xie, Longjie},
title = {Poisson Equation on Wasserstein Space and Diffusion Approximations for Multiscale McKean–Vlasov Equation},
journal = {SIAM Journal on Mathematical Analysis},
volume = {56},
number = {2},
pages = {1495-1524},
year = {2024},
doi = {10.1137/22M1536856},

URL = { 
    
        https://doi.org/10.1137/22M1536856
    
    

},
eprint = { 
    
        https://doi.org/10.1137/22M1536856
    
    

}
}

@article{RSX21,
  author  = {R{\"o}ckner, Michael and Sun, Xiaobin and Xie, Yingchao},
  title   = {Strong Convergence Order for Slow--Fast {McKean--Vlasov}
             Stochastic Differential Equations},
  journal = {Annales de l'Institut Henri Poincar{\'e},
             Probabilit{\'e}s et Statistiques},
  year    = {2021},
  volume  = {57},
  number  = {1},
  pages   = {547--576},
  doi     = {10.1214/20-AIHP1087}
}

@article{WHY24,
  author  = {Wu, Hao and Hu, Junhao and Yuan, Chenggui},
  title   = {Large Deviation for Slow--Fast {McKean--Vlasov}
             Stochastic Differential Equations Driven by Fractional
             Brownian Motions and Brownian Motions},
  journal = {Stochastics and Dynamics},
  year    = {2024},
  volume  = {24},
  number  = {6},
  pages   = {2450044},
  doi     = {10.1142/S0219493724500448}
}

@article{JZ25,
  author  = {Jiang, Yaru and Zhao, Meiling},
  title   = {Averaging Principle for the Fast--Slow {McKean--Vlasov}
             Stochastic Differential Equations Driven by Mixed
             Fractional Brownian Motion},
  journal = {Journal of Applied Probability},
  year    = {2025},
  volume  = {62},
  number  = {4},
  pages   = {1260--1279},
  doi     = {10.1017/jpr.2025.11}
}

@article{HLL22,
  author  = {Hong, Wei and Li, Shihu and Liu, Wei},
  title   = {Strong Convergence Rates in Averaging Principle for
             Slow--Fast {McKean--Vlasov} {SPDEs}},
  journal = {Journal of Differential Equations},
  year    = {2022},
  volume  = {316},
  pages   = {94--135},
  doi     = {10.1016/j.jde.2022.01.021}
}

@article{BS23,
  author  = {Bezemek, Zachary William and Spiliopoulos, Konstantinos},
  title   = {Rate of Homogenization for Fully-Coupled
             {McKean--Vlasov} {SDEs}},
  journal = {Stochastics and Dynamics},
  year    = {2023},
  volume  = {23},
  number  = {2},
  pages   = {2350013},
  doi     = {10.1142/S0219493723500132}
}

@article{SLLG25,
  author  = {Shi, Yangyang and Liu, Hui and Gao, Hongjun},
  title   = {Weak and Strong Averaging Principle for Non-Autonomous
             Slow--Fast {McKean--Vlasov} {SDEs} With Almost Periodic
             Coefficients},
  journal = {Mathematical Methods in the Applied Sciences},
  year    = {2025},
  volume  = {48},
  number  = {16},
  pages   = {15111--15131},
  doi     = {10.1002/mma.70003}
}

@article{HLL26,
  author  = {Hong, Wei and Li, Ge and Li, Shihu},
  title   = {Multi-Scale {McKean--Vlasov} {SDEs}: Moderate Deviation
             Principle in Different Regimes},
  journal = {Annales de l'Institut Henri Poincar{\'e},
             Probabilit{\'e}s et Statistiques},
  year    = {2026},
  volume  = {62},
  number  = {2},
  pages   = {1112--1163},
  doi     = {10.1214/24-AIHP1538}
}

@article{duong2025coarse,
  title={Coarse graining of stochastic differential equations: averaging and projection method},
  author={Duong, Manh Hong and Hartmann, Carsten and Ottobre, Michela},
  journal={arXiv preprint arXiv:2506.14939},
  year={2025}
}

@article{pavliotis2025linearization,
  title={Linearization of ergodic McKean SDEs and applications},
  author={Pavliotis, Grigorios A and Zanoni, Andrea},
  journal={Nonlinearity},
  volume={38},
  number={8},
  pages={085008},
  year={2025},
  publisher={IOP Publishing}
}

@article{K21,
author="Olav, Kallenberg",
title="Foundations of Modern Probability",
journal="Probability Theory and Stochastic Modelling",
ISSN="2199-3130",
publisher="Springer International Publishing",
year="2021",
DOI="10.1007/978-3-030-61871-1",
URL="https://cir.nii.ac.jp/crid/1360021396131395968"
}
\end{small}

\end{document}